\documentclass[a5paper, 10pt]{article}

\usepackage{geometry}
\usepackage{cite, amsmath, amssymb, booktabs, lastpage, graphicx, float, ,authblk}
\usepackage[hidelinks]{hyperref} 
\usepackage{setspace}
\usepackage{amsthm}

\theoremstyle{plain}

\theoremstyle{remark}

\theoremstyle{definition}

\usepackage{amsthm,verbatim}
\usepackage{chemformula}
\usepackage[version=3]{mhchem} %
\usepackage{physics}
\usepackage{xcolor,soul}
\newcommand{\cbox}[2]{\colorbox{#1}{$#2$}}

\newtheorem{thm}{Theorem}[section]
\newtheorem{dfn}[thm]{Definition}
\newtheorem{pro}[thm]{Problem}
\newtheorem{prop}[thm]{Proposition}

\newtheorem{lem}[thm]{Lemma}
\newtheorem{exa}[thm]{Example}

\newcommand{\R}{\mathbb R}

\newcommand{\be}{\beta}

\newcommand{\Ga}{\Gamma}
\newcommand{\de}{\delta}

\newcommand{\la}{\lambda}

\newcommand{\si}{\sigma}
\newcommand{\Si}{\Sigma}
\newcommand{\ti}{\tilde}
\newcommand{\ep}{\varepsilon}

\newcommand{\Vol}{\mathrm{vol}}
\newcommand{\sign}{\mathrm{sign}}

\newcommand{\SO}{\mathrm{SO}}
\newcommand{\SE}{\mathrm{SE}}
\newcommand{\Eu}{\mathrm{E}}
\newcommand{\Or}{\mathrm{O}}

\newcommand{\CIS}{\mathrm{CIS}} 
\newcommand{\CRS}{\mathrm{CRS}} 

\newcommand{\PDD}{\mathrm{PDD}}

\newcommand{\sort}{\mathrm{sort}}
\newcommand{\ORD}{\mathrm{ORD}}
\newcommand{\mORD}{\overline{\ORD}}
\newcommand{\OSD}{\mathrm{OSD}}
\newcommand{\mOSD}{\overline{\OSD}}
\newcommand{\OCD}{\mathrm{OCD}}

\newcommand{\SCD}{\mathrm{SCD}}
\newcommand{\mSCD}{\overline{\SCD}} 
\newcommand{\ASV}{\mathrm{ASV}}
\newcommand{\ASD}{\mathrm{ASD}}
\newcommand{\ASM}{\mathrm{ASM}} 
\newcommand{\CDF}{\mathrm{CDF}}
\newcommand{\ACV}{\mathrm{ACV}}
\newcommand{\ACD}{\mathrm{ACD}}
\newcommand{\ACM}{\mathrm{ACM}}

\newcommand{\EMD}{\mathrm{EMD}}

\newcommand{\BMD}{\mathrm{BMD}}
\newcommand{\BD}{\mathrm{BD}}
\newcommand{\HD}{\mathrm{HD}}

\newcommand{\SPD}{\mathrm{SPD}}

\newcommand{\aff}{\mathrm{aff}}

\newcommand{\bs}{\hfill $\blacksquare$}

\newcommand{\lra}{\leftrightarrow}

\newcommand{\myskip}{\medskip}
\newcommand{\ar}[1]{\overrightarrow{\vb{#1}}}

\newcommand{\colthree}[3]{ \left( \begin{array}{c} 
 #1 \\ #2 \\ #3 \end{array} \right)}
\newcommand{\dettwo}[4]{ \left| \begin{array}{cc} 
 #1 & #2 \\ #3 & #4 \end{array} \right|}

\newcommand{\mathree}[6]{ \left( \begin{array}{cc} 
 #1 & #4 \\ #2 & #5 \\ #3 & #6 \end{array} \right)}
 
\begin{document}

\title{Complete invariants of atomic clouds under rigid motion with Lipschitz continuous metrics in a polynomial time}
\author[1]{Vitaliy Kurlin$^a$} 
\affil[1]{$^a$School of Computer Science and Informatics, University of Liverpool, Liverpool L69 3BX, UK, vitaliy.kurlin@gmail.com} 
 
\date{\today}

\maketitle

\begin{abstract}
A basic representation of any real molecule is a finite cloud of unordered atoms, many of which are chemically indistinguishable.
A natural equivalence on point clouds in any metric space is defined by isometries that are distance-preserving transformations. 
\medskip

In a Euclidean space, any isometry is a composition of translations, rotations, and reflections.
If points are ordered, the isometry class of this cloud is uniquely determined by the matrix of all pairwise distances.
If $m$ points are unordered, a naive metric based on distance matrices needs exponentially many $m!$ permutations.
\medskip

We define a complete invariant for $n$-dimensional clouds of $m$ unordered points under rigid motion, which distinguishes all mirror images in $\R^n$.
The key challenge was to design a distance on invariant values that is Lipschitz continuous under noise and computable in a polynomial time of cloud sizes, for a fixed dimension $n$.
\end{abstract}

\onehalfspacing

\section{Introduction: clouds of unordered points}
\label{sec:intro}

\begin{figure}[h] 
\includegraphics[height=18mm]{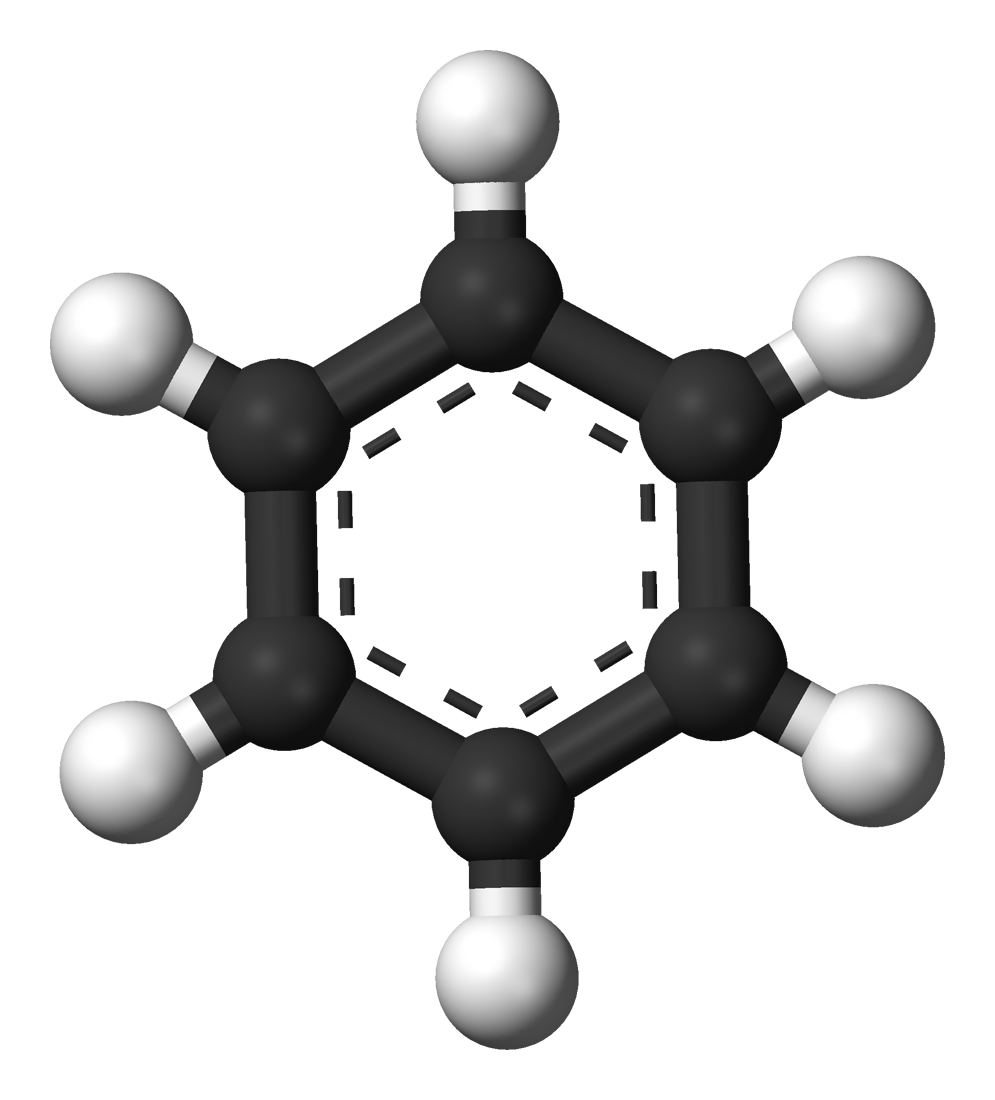}
\hspace*{0.5mm}
\includegraphics[height=18mm]{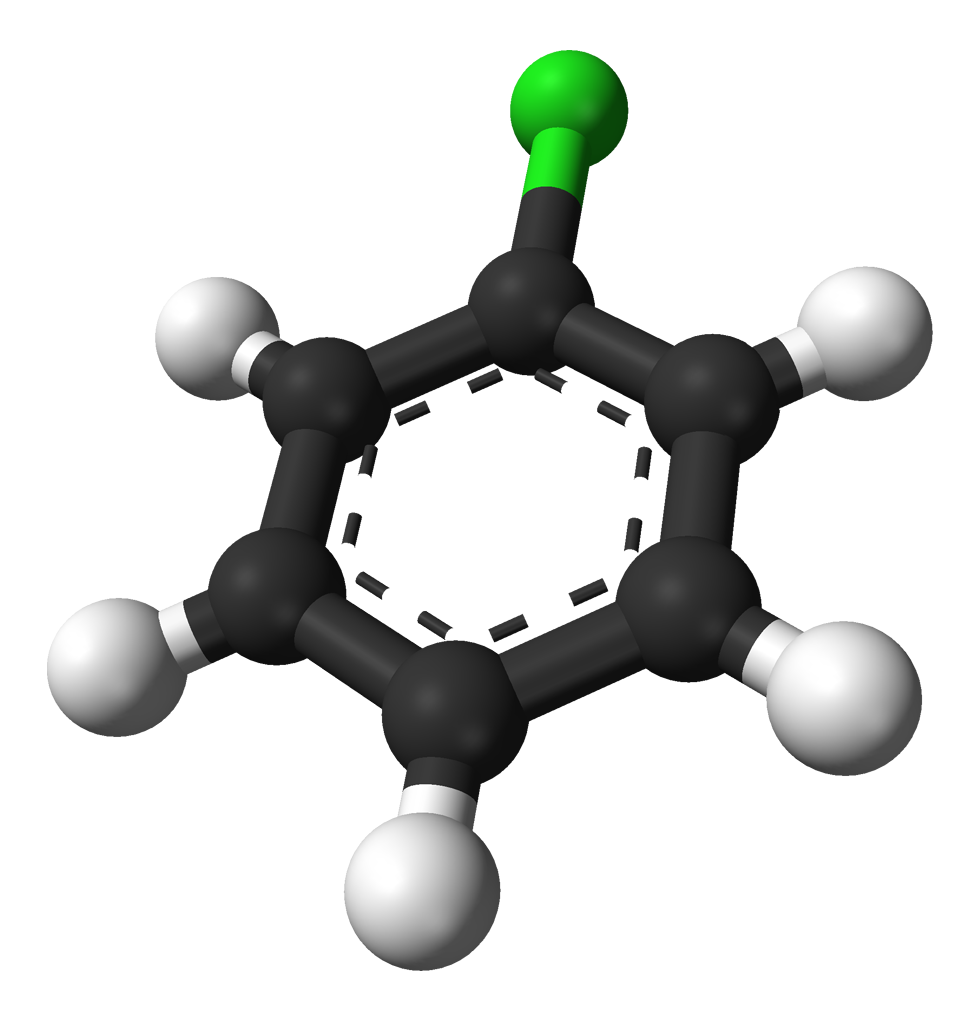}
\hspace*{0.5mm}
\includegraphics[height=18mm]{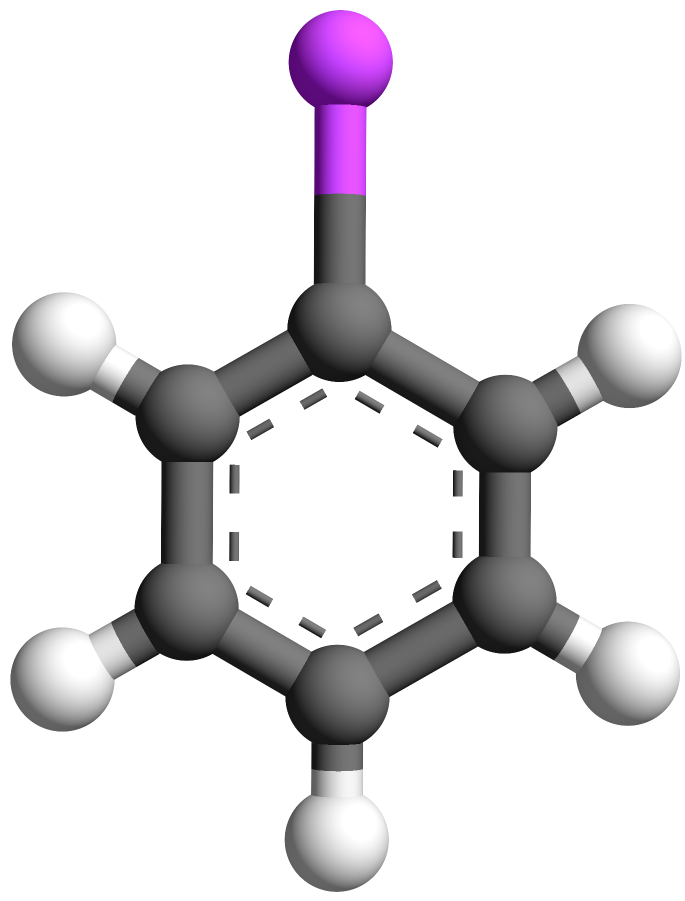}
\hspace*{1mm}
\includegraphics[height=18mm]{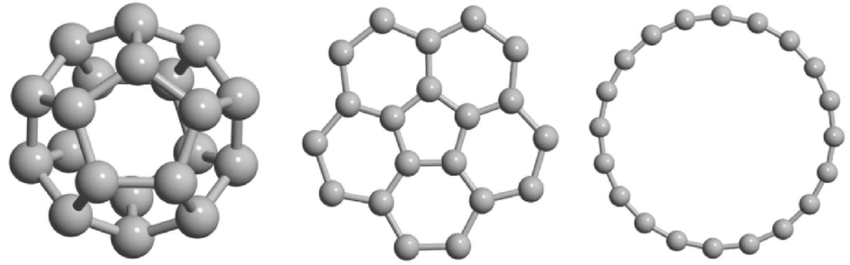}
\caption{
Benzene \ce{C6H6}, phenyllithium \ce{C6 H5 Li}, chlorobenzene \ce{C6 H5 Li}, isomers of \ce{C20} have many indistinguishable atoms.}
\label{fig:molecules}
\end{figure}

Many objects are given as collections of unordered points.
For example, a molecule in $\R^3$ is often represented as a list of points at atomic centers without a canonical order. 
For example, the benzene molecule \ce{C6H6} in the first picture of Fig.~\ref{fig:molecules} consists of 6 carbon atoms and 6 hydrogen atoms, which can be permuted in $(6!)^2>500,000$ ways. 
In computer vision, the ambient dimension is still small, but the size of point clouds can be much larger.
Since points are often unordered in many applications, we study point \emph{clouds} defined as finite sets of unordered (unlabelled) points in $\R^n$.
\medskip

Molecules are often considered as graphs embedded in $\R^3$ with atoms connected by covalent bonds.
However, these bonds are not real sticks in space, and only abstractly represent interatomic interactions that have no precise definitions but depend on many thresholds for distances and angles.
Since atomic nuclei are real physical objects, a cloud of atomic centers remains the most fundamental representation of a real molecule.
\medskip

In addition to the lack of point ordering, another practical challenge is the ambiguity of point representations depending on a coordinate system.
Indeed, any translation or rotation changes a coordinate representation but not the rigid shape of a given cloud.
Clouds $C,C'\subset\R^n$ are called \emph{rigidly equivalent} (denoted by $C\cong C'$) if there is a \emph{rigid motion} (a composition of translations and rotations) from the special Euclidean group $\SE(\R^n)$ that maps $C$ to $C'$.
Rigid motion defines an equivalence relation on clouds.
\medskip

An \emph{equivalence} is any binary relation (denoted by $C\sim C'$) that satisfies the following axioms: (1) \emph{reflexivity:} $C\sim C$, (2) \emph{symmetry:} if $C\sim C'$, then $C'\sim C$, (3) \emph{transitivity:} if $C\sim C'\sim C''$, then $C\sim C''$ for any given clouds.
The transitivity axiom justifies any classification as a splitting into disjoint \emph{equivalence classes}.
Any cloud $C$ generates its class $[C]=\{C'\mid C'\sim C\}$ consisting of all $C'$ equivalent to $C$.
We define the \emph{rigid shape} of any cloud $C\subset\R^n$ as its equivalence class under rigid motion.
\medskip

Another equivalence on point clouds is defined by \emph{isometries}, which are distance-preserving maps and form the full Euclidean group $\Eu(\R^n)$.
The equivalence under isometry (denoted by $C\simeq C'$) is slightly weaker than under rigid motion because any two mirror images belong to the same isometry class but may not be rigidly equivalent.
Compositions of rigid motion or isometry with uniform scaling generate larger groups of transformations and hence smaller spaces of equivalence classes of $m$-point clouds.
These spaces of classes are still high-dimensional and continuous because coordinates of given points, such as atomic centers, are real.
\medskip
 
Rigid shapes can be distinguished by \emph{rigid invariants} $I$ that are descriptors preserved under any rigid motion, i.e. if $C\cong C'$, then $I(C)=I(C')$.
Equivalently, if $I(C)\neq I(C')$, then $C\not\cong C'$ are not rigidly equivalent.
In other words, $I$ has \emph{no false negatives} that are pairs of $C\cong C'$ with $I(C)\neq I(C')$.
For any cloud $C\subset\R^n$, its number $m$ of points is invariant, but the center of mass $O(C)=\dfrac{1}{m}\sum\limits_{p\in C}\vb*{p}$ is not invariant.
Here $\vb*{p}\in\R^n$ denotes the position vector from the origin $0\in\R^n$ to any point $p\in\R^n$.
\medskip

The number of points does not distinguish any non-equivalent clouds of the same size.
A rigid invariant $I$ is \emph{complete} if the converse implication holds: any non-equivalent clouds $C\not\cong C'$ have $I(C)\neq I(C')$.
In other words, $I$ has \emph{no false positives} that are pairs $C\not\cong C'$ with $I(C)=I(C')$. 
\medskip

Recall that the \emph{affine dimension} $\aff(C)$ of a cloud $C$ of points $p_1,\dots,p_m$ in $\R^n$ is the maximum dimension of the vector space generated by all inter-point vectors $\vb*{p_i}-\vb*{p_j}$, $i,j\in\{1,\dots,m\}$.
Then $\aff(C)$ is an isometry invariant and is independent of an order of points of $C$.
Any cloud $C\subset\R^n$ of 2 distinct points has $\aff(C)=1$.
Any cloud $C\subset\R^n$ of 3 points that are not in the same straight line has $\aff(C)=2$.
In the sequel, we assume that any cloud $C\subset\R^n$ is $n$-\emph{dimensional} in the sense that $\aff(C)=n$.
\medskip

A complete invariant answers the first key question (\emph{same or different?}) for objects under a given equivalence.
Since atomic coordinates are always uncertain due to thermal vibrations or experimental noise, the next practical question (\emph{if different, by how much?}) requires the concept of a distance metric \cite{burago2001course} and continuity, which we formalize in Problem~\ref{pro:cloud_invariants_SCD}.

\begin{pro}[complete and continuous invariants of unordered clouds]
\label{pro:cloud_invariants_SCD}
Design an invariant $I$ on all $n$-dimensional clouds of $m$ unordered points in $\R^n$ with values in a metric space to satisfy the following conditions.
\medskip

\noindent
\textbf{(a)}
\emph{Completeness:} 
any $n$-dimensional clouds $C,C'\subset\R^n$ are related by rigid motion (denoted by $C\cong C'$) if and only if $I(C)=I(C')$. 
\medskip

\noindent
\textbf{(b)}
\emph{Reconstruction:} any $n$-dimensional cloud $C\subset\R^n$ of unordered points can be reconstructed from $I(C)$, uniquely under rigid motion.
\medskip

\noindent
\textbf{(c)}
\emph{Metric:} there is a metric $d$ on invariant values satisfying all axioms:
\smallskip

\noindent
(1) $d(I(C),I(C'))=0$ if and only if $C\cong C'$ are related by rigid motion;
\smallskip

\noindent 
(2) \emph{symmetry} : $d(I(C),I(C'))=d(I(C),I(C'))$ for any clouds $C,C'$;
\smallskip

\noindent 
(3) $d(I(C),I(C'))+d(I(C'),I(C''))\geq d(I(C),I(C''))$ for any  $C,C',C''$.
\medskip

\noindent
\textbf{(d)}
\emph{Lipschitz continuity:}
there is a constant $\la>0$ such that, for any $\ep>0$, if a cloud $C'$ is obtained by perturbing every point of a cloud $C\subset\R^n$ up to Euclidean distance $\ep$, then $d(I(C),I(C'))\leq\la\ep$. 
\medskip

\noindent
\textbf{(e)}
\emph{Computability:} for a fixed dimension $n$, the invariant $I(C)$, metric $d(I(C),I(C'))$, and a reconstruction of a cloud $C\subset\R^n$ from $I(C)$ can be found in a polynomial time of the maximum size of any given clouds $C,C'$.
\end{pro}

A reconstruction in~\ref{pro:cloud_invariants_SCD}(b) is stronger than completeness in \ref{pro:cloud_invariants_SCD}(a) because a complete invariant $I(C)$ can be too abstract and may not allow one to reconstruct $C\subset\R^n$. 
In biology, a genetic code is considered complete in practice but does not allow yet a reconstruction of a living organism.
\medskip

The metric axioms in \ref{pro:cloud_invariants_SCD}(c) are important because they are basic requirements for proofs in metric geometry \cite{dorst2010geometric}.
If the first axiom fails, even the zero distance $d\equiv 0$ satisfies the other two axioms.
If the triangle axiom fails or approximately holds with any small additive error, then the $k$-means and DBSCAN clustering are open to adversarial attacks \cite{rass2022metricizing}. 
\medskip

The first metric axiom in \ref{pro:cloud_invariants_SCD}(c) implies the completeness of $I$ in \ref{pro:cloud_invariants_SCD}(a), so the equality $I(C)=I(C')$ is can be checked by verifying if $d(I(C),I(C'))=0$.
For any complete invariant $I$, one can define the discrete metric that takes the value of (say) $1$ for any clouds $C\not\cong C'$, which unhelpfully treats all different clouds (even near-duplicates) as equally distant.
\medskip

The Lipschitz continuity in~\ref{pro:cloud_invariants_SCD}(d) is much stronger than the classical $\ep-\de$ continuity because the Lipschitz constant $\la$ should be independent of a cloud $C$ and a threshold $\ep$.
Indeed, $\dfrac{1}{x}$ is theoretically continuous for all $x>0$, but its behaviour for $x\to 0$ is an example of explosive growth.
\medskip

The computability in~\ref{pro:cloud_invariants_SCD}(e) prevents brute-force attempts to define a complete but rather impractical invariant $I(C)$ by taking the infinite set of images of $C$ under all possible rigid motions or by minimizing a distance $d$ between infinitely many alignments of given clouds $C,C'\subset\R^n$.
All computational complexities are considered in the Random Access Memory (RAM) model, where any value can be accessed in a constant time.
\medskip

The key contribution is the new invariant that completely solves exponentially hard Problem~\ref{pro:cloud_invariants_SCD} for all $n$-dimensional cloud of unordered points.

\section{Past work on invariants of point clouds}
\label{sec:past}
 
This section reviews related approaches with theoretical guarantees.
\medskip

The case of ordered points was resolved more than 90 years ago by reconstructing any sequence $p_1,\dots,p_m\in\R^n$ (uniquely under Euclidean isometry) from the matrix of Euclidean distances $|p_i-p_j|$ or, equivalently, the Gram matrix of scalar products $\vb*{p_i}\cdot\vb*{p_j}$.
Then the difference between $m\times m$ matrices of distances (or ) can be converted into a continuous metric by taking a matrix norm.
In this ordered case, Kendall's shape theory \cite{kendall2009shape} described shape spaces $\Si_m^n$ of $m$-point sequences under isometry in $\R^n$.
\medskip

In the unordered case, comparing $m\times m$ matrices needs an exponential number $m!$ of permutations acting on $m$ given points.
This action by $m!$ permutations makes the invariants and metrics $m!$ times more computationally challenging in Problem~\ref{pro:cloud_invariants_SCD} than for ordered sets.  
\medskip

Cloud invariants based on principal directions were proved to be continuous in general position \cite{kurlin2024polynomial} and inspired the fully separating invariants \cite{hordan2023complete} but degenerate cases, such as 3 points in a line, can lead to discontinuities.  
\medskip
 
Computational geometry studied many metrics between fixed point clouds, which can be theoretically minimized under isometry \cite{huttenlocher1993comparing,chew1999geometric}. 

\begin{dfn}[Minkowski, Hausdorff and bottleneck distances]
\label{dfn:Minkowski+Hausdorff+bottleneck}
\textbf{(a)}
For any real parameter $q\in[1,+\infty)$ and vectors $\vb{u}=(u_1,\dots,u_n)$ and $\vb{v}=(v_1,\dots,v_n)$ in $\R^n$, the \emph{Minkowski} distance is $L_q(\vb{u},\vb{v})=\left(\sum\limits_{i=1}^n |u_i-v_i|^q \right)^{1/q}$.
In the limit case $q=+\infty$, set $L_\infty(\vb{u},\vb{v})=\max\limits_{i=1,\dots,n} |u_i-v_i|$.
\medskip

\noindent
\textbf{(b)}
Let $A,B$ be sets in a space with a metric $d$.
The distance from $p\in A$ to $B$ is $d(p,B)=\inf\limits_{q\in B} d(p,q)$.
Set $d_H(A,B)=\sup\limits_{p\in A} d(p,B)$.
The \emph{Hausdorff} distance is $\HD(A,B)=\max\{d_H(A,B),d_H(B,A)\}$.
The \emph{bottleneck distance} 
$\BD(A,B)=\inf\limits_{g:A\to B} \sup\limits_{p\in A}d(g(p),p)$ is minimized for all bijections $g:A\to B$.
If no bijections $g:A\to B$ exist, set $\BD(A,B)=+\infty$.
\bs
\end{dfn}

The Hausdorff distance between clouds of up to $m$ points, minimized under rigid motion in $\R^2$, is computable in time $O(m^5\log m)$ \cite{chew1997geometric}.
Approximate methods for matching clouds, such as Scale Invariant Feature Transform, sample rotations from the infinite group $\SO(\R^n)$ or use a basis \cite{toews2013efficient}, 
which can be unstable for equal inter-point distances.
The Gromov-Hausdorff distances for 
metric-measure spaces \cite{memoli2011gromov} are defined via infinitely many correspondences between points, but cannot be approximated with a factor less than 3 in polynomial time unless P=NP \cite[Corollary~3.8]{schmiedl2017computational}, see polynomial-time algorithms in partial cases \cite{memoli2021gromov,majhi2024approximating}.
\medskip

Persistent homology \cite{edelsbrunner2022computational} can be considered an isometry invariant of point clouds for filtrations of Vietoris-Rips, Cech, or Delaunay complexes \cite{bauer2014morse}.
The persistence is Lipschitz continuous \cite{cohen2005stability} in the bottleneck distance, as in condition~\ref{pro:cloud_invariants_SCD}(c), but its fastest versions in dimensions 0 and 1 do not distinguish generic families of non-isometric clouds  \cite{smith2024generic}.
Comparing such infinite-size invariants requires sampling the infinite group $\Or(\R^n)$.
\medskip

Problem~\ref{pro:cloud_invariants_SCD} is fully solved for $m=3$ points by the side-side-side theorem saying that any triangles are congruent (isometric) if and only if they have the same triple of sides.
In 2004, Boutin and Kemper \cite{boutin2004reconstructing} proved in the unordered case that all sorted pairwise distances in Definition~\ref{dfn:SPD} uniquely determine almost any (generic) cloud $C\subset\R^n$ under isometry.

\begin{dfn}[Sorted Pairwise Distances $\SPD$]
\label{dfn:SPD}
For any cloud $C$ of $m$ unordered points, the vector $\SPD(C)$ of \emph{Sorted Pairwise Distances} has $\dfrac{m(m-1)}{2}$ distances between points of $C$, written in increasing order.
\bs
\end{dfn}

\begin{figure}[h!]
\centering
\includegraphics[width=\linewidth]{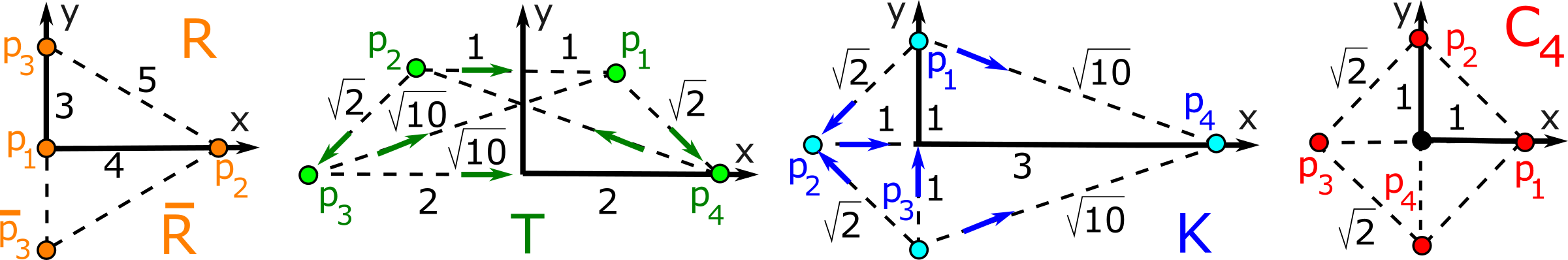}
\caption{\textbf{1st}: the right-angled cloud $R\subset\R^2$ consists of points $p_1=(0,0)$, $p_2=(4,0)$, $p_3=(0,3)$, and its mirror image $\bar R$ of $p_1,p_2$, and $\bar p_3=(0,-3)$ with respect to the $x$-axis. 
\textbf{2nd}: the trapezoid cloud $T\subset\R^2$ consists of points $p_1=(1,1)$, $p_2=(-1,1)$, $p_3=(-2,0)$, $p_4=(2,0)$. 
\textbf{3rd}: the kite cloud $K\subset\R^2$ consists of points $p_1=(0,1)$, $p_2=(-1,0)$, $p_3=(0,-1)$, $p_4=(3,0)$. 
\textbf{4th}: the square cloud $S\subset\R^2$ consists of 
$p_1=(1,0)$, $p_2=(0,1)$, $p_3=(-1,0)$, $p_4=(0,-1)$. 
}
\label{fig:clouds_oriented}
\end{figure}

The invariant $\SPD(C)$ is incomplete already for $m=4$ points in $\R^2$, see the non-isometric 4-point clouds $T\not\cong K$ in the two middle pictures of Fig.~\ref{fig:clouds_oriented}, which have $\SPD(T)=[\sqrt{2},\sqrt{2},2,\sqrt{10},\sqrt{10},4]=\SPD(K)$.

Problem~\ref{pro:cloud_invariants_SCD} was solved for $m=4$ unordered points under isometry by the stronger invariant $\PDD$ in Definition~\ref{dfn:PDD_finite} by \cite[Theorem~5.3]{widdowson2026pointwise}, but the case of $m\geq 5$ unordered points remained open even in the plane $\R^2$.

\begin{dfn}[Pointwise Distance Distribution $\PDD(C;k)$]
\label{dfn:PDD_finite}
In a metric space, let $C$ be any finite cloud of unordered points $p_1,\dots,p_m$, which are labeled here only for convenience.
Fix an integer $1\leq k<m$.
For every point $p_i\in C$, let $d_1(p_i)\leq\dots\leq d_k(p_i)$ be ordered distances from $p_i$ to its $k$ nearest neighbors in $C$.
The matrix $D(C;k)$ has $m$ rows consisting of the distances $d_1(p_i),\dots,d_k(p_i)$ for $i=1,\dots,m$.
If any $l\geq 2$ rows coincide, we collapse them into a single row with the weight $\dfrac{l}{m}$.
The resulting distribution of unordered rows, each having a weight and $k$ ordered distances, is called the \emph{Pointwise Distance Distribution} $\PDD(C;k)$ \cite{widdowson2022resolving}. 
\bs
\end{dfn}

If we consider only point clouds in a general position away from singular configurations, one can order all points by their distinct distances to the center of mass.
In this general position, Problem~\ref{pro:cloud_invariants_SCD} reduces to the easier ordered case but remains exponentially hard for degenerate clouds. 
\medskip

\begin{figure}[h!]
\includegraphics[width=\textwidth]{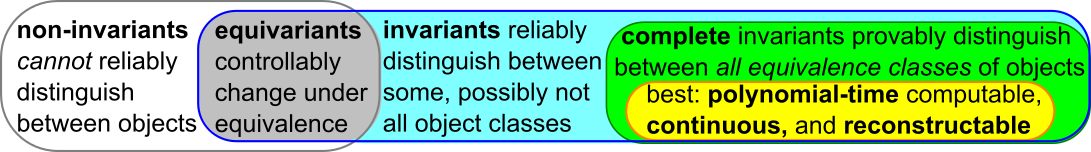}
\caption{
Non-invariant descriptors, such as point coordinates, are unreliable for comparisons under isometry.
Any isometry $f$ changes equivariant descriptors, such as the center of mass, by the same $f$.
Only invariants preserved under all equivalences can classify objects.
New invariants allow reconstructions and polynomial-time Lipschitz continuous metrics. 
}
\label{fig:invariants}
\end{figure}  

Geometric Deep Learning \cite{bronstein2021geometric} experimentally tries to solve Problem~\ref{pro:cloud_invariants_SCD} by neural networks \cite{kondor2018generalization} that output invariants or equivariants. 
An \emph{equivariant} descriptor $E$ satisfies $E(f(A))=T_f(E(A))$, where $T_f$ depends only on a transformation $f$, not on a cloud $A\subset\R^n$, but $T_f$ may not be the identity, as needed for \emph{invariants}, to distinguish clouds under equivalences \cite{satorras2021n}. 
More hidden parameters usually improve accuracy at a higher cost \cite{dong2018boosting}, 
but still do not guarantee completeness, which requires \emph{no false negatives} and \emph{no false positives} for all possible data. 
The universal approximation property of neural networks 
means the completeness of infinite-size invariants \cite{maron2019universality}, 
while Problem~\ref{pro:cloud_invariants_SCD} requires finite polynomial-time invariants.

\section{Oriented Simplexwise Distribution ($\OSD$)}
\label{sec:OSD}

This section starts solving Problem~\ref{pro:cloud_invariants_SCD} by
introducing the Oriented Simplexwise Distribution ($\OSD$) in Definition~\ref{dfn:OSD} after a few auxiliary concepts below.
For any ordered vectors $\vb*{v_1},\dots,\vb*{v_n}\in\R^n$, let $\sign(\vb*{v_1},\dots,\vb*{v_n})$ denote the sign of the $n\times n$ determinant with the column vectors $\vb*{v_1},\dots,\vb*{v_n}$, which can be $\pm 1$ or $0$ (if the vectors $\vb*{v_1},\dots,\vb*{v_n}$ are linearly dependent).

\begin{dfn}[matrix $D(A)$ and cloud $M(C;A)\subset\R^{n+1}$ for $A\subseteq C$]
\label{dfn:D(A)+M(C;A)}
Let $C$ be any $n$-dimensional cloud of $m$ unordered points in $\R^n$.
For a \emph{base sequence} $A$ of $n$ ordered points $p_1,\dots,p_n\in C$, let $D(A)$ be the $n\times n$ matrix of distances $|\vb*{p_i}-\vb*{p_j}|$ for $i,j=1,\dots,n$.
For any point $q\in C\setminus A$, set $s(q;A)=\sign(\vb*{q}-\vb*{p_1},\dots,\vb*{q}-\vb*{p_n})$.
Let $d(q)\in\R^{n+1}$ be the point with the coordinates $|\vb*{q}-\vb*{p_1}|,\dots,|\vb*{q}-\vb*{p_n}|,s(q;A)$. 
Let $M(C;A)\subset\R^{n+1}$ be the cloud of $m-n$ unordered points $d(q)$ for all points $q\in C\setminus A$.
\bs
\end{dfn}

Since all points $d(q)\in M(C;A)\subset\R^{n+1}$ have fixed coordinates, it will be convenient to represent the cloud $M(C;A)$ by a matrix with $m-n$ permutable columns and $n+1$ ordered rows.
Let $S_n$ denote the symmetric group on the indices $1,\dots,n$.
Any permutation $\xi\in S_n$ is a composition of some $t$ transpositions $i\lra j$ and has $\sign(\xi)=(-1)^t$. 

\begin{dfn}[oriented distributions $\ORD(C;A)$ and $\OSD(C)$ for a cloud $C\subset\R^n$]
\label{dfn:OSD}
In the notations of Definition~\ref{dfn:D(A)+M(C;A)}, any permutation $\xi\in S_n$ acts on the matrix $D(A)$ by mapping $D(A)_{i,j}$ to $D(A)_{\xi(i),\xi(j)}$ and on the cloud $M(C;A)$ of points $d(q)$ for $q\in C\setminus A$ by permuting the first $n$ coordinates of $d(q)$ and multiplying the $(n+1)$-st coordinate $s(q;A)$ by $\sign(\xi)$.
The \emph{Oriented Relative Distribution} $\ORD(C;A)$ is the equivalence class of the pair $[D(A);M(C;A)]$ under all permutations $\xi\in S_n$, which act on both $D(A)$ and $M(C;A)$. 
The \emph{Oriented Simplexwise Distribution} $\OSD(C)$ is the unordered collection of $\ORD(C;A)$ for all $\dbinom{m}{n}$ unordered subsets $A\subseteq C$ of $n$ points. 
The \emph{mirror images} $\mORD(C;A)$ and $\mOSD(C)$ are obtained from $\ORD(C;A)$ and $\OSD(C)$, respectively, by reversing the signs of $(n+1)$-st coordinates of points $d(q)\in M(C;A)$ for $q\in C\setminus A$.
\bs
\end{dfn}

\begin{exa}[$\OSD$ for clouds in Fig.~\ref{fig:clouds_oriented}]
\label{exa:OSD}
\textbf{(a)}
In $\R^2$ with the counter-clockwise orientation, the right-angled cloud $R$ consisting of the vertices $p_1=(0,0)$, $p_2=(4,0)$, $p_3=(0,3)$ of the triangle in Fig.~\ref{fig:clouds_oriented}~(1st) has the distribution $\OSD(R)$ consisting of 
$\ORD(R;(p_1,p_2))=\left[4,\colthree{3}{5}{+1}\right],$
$\ORD(R;(p_2,p_3))=\left[5,\colthree{4}{3}{+1}\right],$
$\ORD(R;(p_3,p_1))=\left[3,\colthree{5}{4}{+1}\right].$  
\smallskip

If we swap the points $p_1\lra p_3$, the last $\ORD$ above changes to the equivalent form $\ORD(R;(p_1,p_3))=\left[3,\colthree{4}{5}{-}\right]$, without affecting others.
If we reflect $R$ across the $x$-axis, the mirror image $\bar R$ of $p_1$, $p_2$, $\bar p_3=(0,-3)$ has 
$\OSD(\bar R)=\mOSD(R)$ with 
$\ORD(\bar R;(p_1,p_2))=\left[4,\colthree{3}{5}{-1}\right],$
$\ORD(\bar R;(p_2,\bar p_3))=\left[5,\colthree{4}{3}{-}\right],$
$\ORD(\bar R;(p'_3,p_1))=\left[3,\colthree{5}{4}{-1}\right],$
 which differs from $\OSD(R)$ even if we swap the points in each pair.
\medskip

\noindent
\textbf{(b)}
Table~\ref{tab:OSD+TK} shows the Oriented Simplexwise Distributions for the clouds $T,K$ in Fig.~\ref{fig:clouds_oriented}.
Each row has similar $\ORD$s with \hl{highlighted} differences.
\end{exa}

\begin{table}[h!]
  \centering
  \begin{tabular}{@{}ll@{}}
    $\ORD$s in $\OSD(\cbox{yellow}{T})$ & 
    $\ORD$s in $\OSD(\cbox{yellow}{K})$ \\
    \myskip
    
	$\left[\sqrt{2},\mathree{2}{\cbox{yellow}{\sqrt{10}}}{-1}{\sqrt{10}}{4}{-1}\right]$ &
	$\left[\sqrt{2},\mathree{2}{\cbox{yellow}{\sqrt{2}}}{-1}{\sqrt{10}}{4}{-1}\right]$ 
\myskip \\

	$\left[\sqrt{2},\mathree{2}{\cbox{yellow}{\sqrt{10}}}{+1}{\sqrt{10}}{4}{+1}\right]$ &
	$\left[\sqrt{2},\mathree{2}{\cbox{yellow}{\sqrt{2}}}{+1}{\sqrt{10}}{4}{+1}\right]$ 
\myskip \\
	
	    $\left[2,\mathree{\sqrt{2}}{\cbox{yellow}{\sqrt{10}}}{-1}{\sqrt{10}}{\cbox{yellow}{\sqrt{2}}}{\cbox{yellow}{-1}}\right]$ & 
    $\left[2,\mathree{\sqrt{2}}{\cbox{yellow}{\sqrt{2}}}{-1}{\sqrt{10}}{\cbox{yellow}{\sqrt{10}}}{\cbox{yellow}{+1}}\right]$
\myskip \\
    
	$\left[\sqrt{10},\mathree{\sqrt{2}}{\cbox{yellow}{2}}{\cbox{yellow}{+1}}{\cbox{yellow}{4}}{\cbox{yellow}{\sqrt{2}}}{-1}\right]$ &
	$\left[\sqrt{10},\mathree{\sqrt{2}}{\cbox{yellow}{4}}{\cbox{yellow}{-1}}{\cbox{yellow}{2}}{\cbox{yellow}{\sqrt{10}}}{-1}\right]$ \myskip \\
	
	$\left[\sqrt{10},\mathree{\sqrt{2}}{\cbox{yellow}{2}}{\cbox{yellow}{-1}}{\cbox{yellow}{4}}{\cbox{yellow}{\sqrt{2}}}{+1}\right]$ & 
	$\left[\sqrt{10},\mathree{\sqrt{2}}{\cbox{yellow}{4}}{\cbox{yellow}{+1}}{\cbox{yellow}{2}}{\cbox{yellow}{\sqrt{10}}}{+1}\right]$ \smallskip \\
	
	    $\left[4,\mathree{\sqrt{2}}{\sqrt{10}}{+1}{\cbox{yellow}{\sqrt{10}}}{\cbox{yellow}{\sqrt{2}}}{\cbox{yellow}{+1}}\right]$ &
	$\left[4,\mathree{\sqrt{2}}{\sqrt{10}}{+1}{\cbox{yellow}{\sqrt{2}}}{\cbox{yellow}{\sqrt{10}}}{\cbox{yellow}{-1}}\right]$ 
\end{tabular}
 \caption{The Oriented Simplexwise Distributions $\OSD$s from Definition~\ref{dfn:OSD} for the 4-point clouds $T,K\subset\R^2$ in Fig.~\ref{fig:clouds_oriented}.
}
  \label{tab:OSD+TK}
\end{table}

Lemma~\ref{lem:distance_realisation} provides a criterion for a matrix to be realizable by squared distances between any $m$ ordered points in $\R^n$.
\smallskip

\begin{lem}[realization of distances]
\label{lem:distance_realisation}
\textbf{(a)}
A symmetric $m\times m$ matrix of $s_{ij}\geq 0$ with $s_{ii}=0$ is realizable as a matrix of squared distances between $p_0=0,p_1,\dots,p_{m-1}\in\R^n$ for some $n$ \emph{if and only if} the $(m-1)\times(m-1)$ matrix $G$ of $g_{ij}=\dfrac{s_{0i}+s_{0j}-s_{ij}}{2}$ has only non-negative eigenvalues.
\smallskip

\noindent
\textbf{(b)}
If $G$ has only non-negative eigenvalues, then $\aff(0,p_1,\dots,p_{m-1})$ equals the number $k\leq m-1\leq n$ of positive eigenvalues of $G$. 
In this case, $g_{ij}=p_i\cdot p_j$ define the \emph{Gram matrix} of the vectors $\vb*{p_1},\dots,\vb*{p_{m-1}}\in\R^n$, which are uniquely determined in time $O(m^3)$ under a map from $\Or(\R^n)$.
\end{lem}
\newcommand{\lemdistancerealisation}{
\begin{proof}
We extend \cite[Theorem~1]{dekster1987edge} for $m\leq n$ and find points $p_1,\dots,p_{m-1}$ in $\R^n$ in time $O(m^3)$, uniquely under an orthogonal map from $\Or(\R^n)$.
\myskip

The part \emph{only if} $\Rightarrow$.
Let a symmetric matrix $S$ consist of squared distances between the points $p_0=0,p_1,\dots,p_{m-1}\in\R^n$.
For $i,j=1,\dots,m-1$, 
$g_{ij}=\dfrac{s_{0i}+s_{0j}-s_{ij}}{2}=\dfrac{|\vb*{p_i}|^2+|\vb*{p_j}|^2-|\vb*{p_i}-\vb*{p_j}|^2}{2}=\vb*{p_i}\cdot\vb*{p_j}$ form the Gram matrix, which can be written as $G=P^T P$, where the columns of the $n\times(m-1)$ matrix $P$ are the vectors $\vb*{p_1},\dots,\vb*{p_{m-1}}$.
For any vector $\vb*{v}\in\R^{m-1}$, we have 
$0\leq |P \vb*{v}|^2=(P\vb*{v})^T(P\vb*{v})=\vb*{v}^T(P^T P)\vb*{v}=\vb*{v}^T G \vb*{v}.$
Since the quadratic form $\vb*{v}^T G \vb*{v}\geq 0$ for any $\vb*{v}\in\R^{m-1}$, the matrix $G$ is positive semi-definite, i.e. $G$ has non-negative eigenvalues \cite[Theorem~7.2.7]{horn2012matrix}.  
\smallskip

The part \emph{if} $\Leftarrow$.
For any positive semi-definite matrix $G$, there is an orthogonal matrix $Q$ such that $Q^T G Q=D$ is the diagonal matrix, whose $m-1$ diagonal elements are non-negative eigenvalues of $G$.
The diagonal matrix $\sqrt{D}$ consists of the square roots of the eigenvalues of $G$.
The number of positive eigenvalues of $G$ equals the dimension $k=\aff(\{0,p_1,\dots,p_{m-1}\})$ of the subspace that is linearly spanned by $\vb*{p_1},\dots,\vb*{p_{m-1}}$. 
\smallskip

We may assume that all $k\leq n$ positive eigenvalues of $G$ correspond to the first $k$ coordinates of $\R^n$.
Since $Q^T=Q^{-1}$, the matrix $G=Q D Q^T=(Q\sqrt{D})(Q\sqrt{D})^T$ becomes the Gram matrix of the columns of $Q\sqrt{D}$.
These columns become the reconstructed points $p_1,\dots,p_{m-1}\in\R^n$.
\smallskip

If there is another diagonalization $\ti Q^T G \ti Q=\ti D$ for $\ti Q\in \Or(\R^n)$, then $\ti D$ differs from $D$ by a permutation of eigenvalues, which is realized by an orthogonal map, so we set $\ti D=D$.
Then $G=\ti Q D \ti Q^T=(\ti Q\sqrt{D})(\ti Q\sqrt{D})^T$ is the Gram matrix of the columns of $\ti Q\sqrt{D}$.
The new columns differ from the previously reconstructed points $p_1,\dots,p_{m-1}\in\R^n$ by the orthogonal map $Q\ti Q^T$.
Hence the reconstruction is unique under orthogonal transformations from $\Or(\R^n)$.
Computing $\vb*{p_1},\dots,\vb*{p_{m-1}}$ as eigenvectors needs a diagonalization of $G$ in time $O(m^3)$, see  \cite[section~11.5]{press2007numerical}.
\end{proof}
}\lemdistancerealisation

Lemma~\ref{lem:seq_reconstruction}(a) holds for all clouds, including degenerate ones, e.g. for 3 points in a straight line.
Any $n-1$ points have $\aff(\{p_1,\dots,p_{n-1}\})\leq n-2$.
For example, any two distinct points in $A\subset\R^3$ generate a straight line. 

\begin{lem}[sequence reconstruction]
\label{lem:seq_reconstruction}
\textbf{(a)}
Any sequence of ordered points $p_1,\dots,p_m$ in $\R^n$ can be reconstructed (uniquely under isometry) from the matrix of the Euclidean distances $|\vb*{p_i}-\vb*{p_j}|$ in time $O(m^3)$.
If all distances are divided by $R=\max\limits_{i=1,\dots,m}|\vb*{p_i}|$, the reconstruction of points $p_1,\dots,p_m$ is unique under isometry and uniform scaling in $\R^n$.
\medskip

\noindent
\textbf{(b)}
If $m\leq n$, the uniqueness of reconstructions in part (a) remains true if we replace isometry with rigid motion in $\R^n$.
\end{lem}
\newcommand{\lemreconstruction}{
\begin{proof}
\textbf{(a)}
By translation, we can fix the point $p_1$ at the origin $0\in\R^n$.
Let $G$ be the $(m-1)\times (m-1)$ matrix $g_{ij}=\dfrac{|\vb*{p_i}|^2+|\vb*{p_j}|^2-|\vb*{p_i}-\vb*{p_j}|^2}{2}=\vb*{p_i}\cdot \vb*{p_j}$, where $i,j=2,\dots,m$, which is obtained from the squared distances between the points $p_1=0,p_2,\dots,p_m$.
\smallskip

By Lemma~\ref{lem:distance_realisation} if $G$ has $k\leq n$ positive eigenvalues, then $p_1=0,\dots,p_m$ can be uniquely determined under isometry in $\R^k\subset\R^n$ in time $O(m^3)$.
If all distances are divided by the same radius $R=\max\limits_{i=1,\dots,m}|\vb*{p_i}|$, the construction guarantees uniqueness under isometry and uniform scaling. 
\medskip

\noindent
\textbf{(b)}
If $m\leq n$, any mirror images of $p_1,\dots,p_m\in\R^n$, after a suitable rigid motion, can be assumed to belong to an $(n-1)$-dimensional hyperspace $H\subset\R^n$, where they are matched by a mirror reflection $H\to H$ with respect to an $(n-2)$-dimensional subspace $S\subset H$, which is realized by the $180^\circ$ orientation-preserving rotation around $S$.
\end{proof}
}\lemreconstruction

Lemma~\ref{lem:seq_reconstruction}(b) for $m=n=3$ implies that any triangle is determined by its sides, uniquely under rigid motion in $\R^3$.
For example, sides $3,4,5$ define a right-angled triangle whose mirror images are not related by rigid motion within a plane $H\subset\R^3$, but are matched by a rigid motion in $H$ and a $180^\circ$ rotation of $\R^3$ around a line in the plane $H$.   
\smallskip

Lemma~\ref{lem:ORD_reconstruction} below will help guarantee the reconstruction in Problem~\ref{pro:cloud_invariants_SCD}.

\begin{lem}[reconstruction from $\ORD$]
\label{lem:ORD_reconstruction}
Any $n$-dimensional cloud $C\subset\R^n$ of $m$ unordered points can be reconstructed, uniquely under rigid motion, from $\ORD(C;A)$ for any base sequence $A\subseteq C$ with $\aff(A)=n-1$.   
\end{lem}
\begin{proof}
By Lemma~\ref{lem:seq_reconstruction}(b), any base sequence $A\subseteq C$ of $n$ points can be reconstructed, uniquely under rigid motion in $\R^n$, from the distance matrix $D(A)$ in Definition~\ref{dfn:OSD}.
We may assume that the $n$ points $p_1,\dots,p_n$ of $A\subseteq C$ span the subspace of the first $n-1$ coordinate axes of $\R^n$. 
\smallskip

We will prove that any point $q\in C\setminus A$ has a unique location in $\R^n$, determined by the $n$ distances $|\vb*{q}-\vb*{p_1}|,\dots,|\vb*{q}-\vb*{p_n}|$ written in a column of the matrix $\ORD(C;A)$.
Since the points of $A$ do not belong to any $(n-1)$-dimensional affine subspace of $\R^n$, the $n$ spheres $S(p_i;|\vb*{q}-\vb*{p_i}|)$ with radii $|\vb*{q_i}-\vb*{p_i}|$ and centers $p_i$, $i=1,\dots,n$, contain $q$ and their full intersections consists of one or two points. 
We can uniquely choose $q$ among these two options due to the sign of the determinant on the column vectors $\vb*{q}-\vb*{p_1},\dots,\vb*{q}-\vb*{p_n}$ in the bottom row of $\ORD(C;A)$. 
\end{proof}

Lemma~\ref{lem:ORD_reconstruction} implies that $\ORD(C;A)$ can have identical columns only for degenerate subsets $A\subseteq C$ with $\aff(A)<n-1$.
For example, let $n=3$ and $A$ consist of three points $p_1,p_2,p_3$ in the same straight line $L\subset\R^3$.
The three distances $|\vb*{q}-\vb*{p_i}|$, $i=1,2,3$, to any other point $q\in C$ outside the line $L$ define three spheres $S(p_i;|\vb*{q}-\vb*{p_i}|)$ that share a common circle in $\R^3$, so the position of $q$ is not uniquely determined in this case.
\smallskip

Though one $\ORD(C;A)$ with $\aff(A)=n-1$ suffices to reconstruct $C\subset\R^n$ uniquely under rigid motion, the dependence on $A$ requires us to consider the Oriented Simplexwise Distribution $\OSD(C)$ of $\ORD(C;A)$ for all $n$-point subsets $A\subseteq C$ to get a complete invariant in Theorem~\ref{thm:OSD_complete}. 
Equality $\OSD(C)=\OSD(C')$ between unordered distributions is interpreted as a bijection $\OSD(C)\to\OSD(C')$ matching all ORDs. 

\begin{thm}[time and completeness of $\OSD$ under rigid motion]
\label{thm:OSD_complete}
\textbf{(a)}
For any $n$-dimensional cloud $C\subset\R^n$ of $m$ unordered points, the Oriented Simplexwise Distribution $\OSD(C)$ can be computed in time $O(m^{n+1})$. 
\medskip

\noindent
\textbf{(b)}
Any $n$-dimensional clouds $C,C'\subset\R^n$ of $m$ unordered points can be exactly matched by rigid motion \emph{if and only if} $\OSD(C)=\OSD(C')$.
\medskip

\noindent
\textbf{(c)}
Any $n$-dimensional clouds $C,C'\subset\R^n$ of $m$ unordered points can be exactly matched by isometry \emph{if and only if} $\OSD(C)$ equals $\OSD(C')$ or its mirror image $\mOSD(C')$.
\end{thm}
\begin{proof}
\textbf{(a)}
To compute $\OSD(C)$, we consider $\dbinom{m}{n}=\dfrac{m!}{n!(m-n)!}$ subsets $A\subseteq C$ of $n$ points.
For each fixed $A$, the matrix $D(A)$ of $\dfrac{n(n-1)}{2}$ pairwise distances needs $O(n^2)$ time.
The Oriented Relative Distribution $\ORD(C;A)=[D(A);M(C;A)]$ includes $n(m-n)$ distances and $m-n$ signs, where each sign requires determinant computations in time $O(n^3)$ by Gaussian elimination.
So $\ORD(C;A)$ is computed in time $O(n^3m)$.
Multiplying the last time by the number $\dbinom{m}{n}=\dfrac{m!}{n!(m-n)!}$ of $n$-point subsets $A\subseteq C$, we estimate the final time for $\OSD(C)$ as
$$O\left(\dfrac{m!n^3m}{n!(m-n)!}\right)=O\left(\dfrac{m^2(m-1)\dots(m-n+1)n^2}{(n-1)!}\right)=O(m^{n+1}).$$

\noindent
\textbf{(b)}
Part \emph{if} $\Leftarrow$.
Any bijection $\OSD(C)\to\OSD(C')$ matches $\ORD(C;A)$ with $\ORD(C';A')$ for some base sequences $A\subseteq C$ and $A'\subseteq C'$. 
Any $n$-dimensional cloud $C\subset\R^n$ has a base sequence $A$ with $\aff(A)=n-1$.
By Lemma~\ref{lem:ORD_reconstruction} an equality $\ORD(C;A)=\ORD(C';A')$ for base sequences $A,A'$ with $\aff=n-1$ guarantees that $C,C'$ are related by rigid motion. 
\smallskip

Part \emph{only if} $\Rightarrow$.
Any rigid motion $f$ of $\R^n$ bijectively maps $C$ to $f(C)=C'$ and 
hence induces a bijection $\OSD(C)\to\OSD(C')$. 
\medskip

\noindent
\textbf{(c)}
follows from part (b).
Indeed, any orientation-reversing isometry $C\to C'$ induces a bijection $\OSD(C)\to\mOSD(C')=\OSD(\bar C')$
Conversely, $\OSD(C)=\mOSD(C')$ implies that $\ORD(C;A)=\ORD(\bar C';A')$ for some base sequences $A\subseteq C$ and $A'\subset\bar C'$, where $\bar C'$ is a mirror image of $C'$.
Then $C,\bar C'$ are related by rigid motion due to Lemma~\ref{lem:ORD_reconstruction}, so $C\simeq C'$.
\end{proof}

For a fixed dimension $n$, the time $O(m^{n+1})$ is polynomial in the number $m$ of points and will be improved to $O(m^n)$ for a smaller invariant below. 

\section{Simplexwise Centered Distribution ($\SCD$)}
\label{sec:SCD}

This section simplifies the $\OSD$ invariant to the Simplexwise Centered Distribution ($\SCD$) in Definition~\ref{dfn:SCD}.
The Euclidean structure of $\R^n$ allows us to translate the \emph{center of mass} $O(C)=\dfrac{1}{m}\sum\limits_{p\in C} \vb*{p}$ of a given $m$-point cloud $C\subset\R^n$ to the origin $0\in\R^n$.
Then Problem~\ref{pro:cloud_invariants_SCD} reduces to the invariance under rotations from the special orthogonal group $\SO(\R^n)$.
\myskip

Definition~\ref{dfn:OSD} introduced the Oriented Simplexwise Distribution $\OSD(C)$ as an ordered collection of $\ORD$s depending on $n$-point subsets of $C$.
Including the center of mass $O(C)$ allows us to consider the smaller number of $\dbinom{m}{n-1}$ subsets $A\subseteq C$ consisting of $n-1$ points instead of $n$ points.
\myskip

Though the center of mass is uniquely determined for any cloud $C\subset\R^n$ of unordered points, real applications may offer one or several labelled points of $C$ that substantially speed up metrics on invariants.
\smallskip

For example, an atomic neighborhood in a solid material is a cloud $C\subset\R^3$ of atoms around a central atom, which may not be the center of mass of $C$, but can be an extra base point in all constructions below.

\begin{table}[h!]
\centering
\begin{tabular}{lll}
$\BD(A,B)$ & Bottleneck Distance between sets $A,B$ & Definition~\ref{dfn:Minkowski+Hausdorff+bottleneck} \\
$\OCD(C;A)$ & Oriented Centered Distribution for $A\subset C$ & Definition~\ref{dfn:SCD} \\
$\ORD(C;A)$ & Oriented Relative Distribution for $A\subset C$ & 
Definition~\ref{dfn:OSD} \\
$\OSD(C)$ & Oriented Simplexwise Distribution of $C$ & 
Definition~\ref{dfn:OSD} \\
$\PDD(C;k)$ & Pointwise Distance Distribution of $C$  & Definition~\ref{dfn:PDD_finite} \\
$\SCD(C)$ & Simplexwise Centered Distribution of $C$ & 
Definition~\ref{dfn:SCD} \\
$\SPD(C)$ & Sorted Pairwise Distances of a cloud $C$  & Definition~\ref{dfn:SPD}
\end{tabular}
\caption{Acronyms and references for concepts in sections~\ref{sec:past}, \ref{sec:OSD}, and \ref{sec:SCD}.}
\label{tab:acronyms_SCD_invariants}
\end{table}

\begin{dfn}[Simplexwise Centered Distribution $\SCD$]
\label{dfn:SCD}
Let $C\subset\R^n$ be any cloud of $m$ unordered points.
For any base sequence $A$ of $n-1$ ordered points $p_1,\dots,p_{n-1}\in C$, set $\ti A=A\cup\{O(C)\}$ by adding the center of mass $O(C)$ as the last $n$-th point.
Following Definition~\ref{dfn:D(A)+M(C;A)}, take the $n\times n$ distance matrix $D(\ti A)$ and the cloud $M(C;\ti A)\subset\R^{n+1}$ of points $d(q)$ for $q\in C\setminus A$.
Any permutation $\xi\in S_{n-1}$ acts on $D(\ti A)$ by mapping $D(\ti A)_{i,j}$ to $D(\ti A)_{\xi(i),\xi(j)}$ for $i,j=1,\dots,n-1$, and on $M(C;\ti A)\subset\R^{n+1}$ by permuting the first $n-1$ coordinates and by multiplying the $(n+1)$-st coordinate by $\sign(\xi)$.
The \emph{Oriented Centered Distribution} $\OCD(C;A)$ is the equivalence class of pairs $[D(\ti A),M(C;\ti A)]$ under all permutations $\xi\in S_{n-1}$ on points of $A$.
The \emph{Simplexwise Centered Distribution} $\SCD(C)$ is the unordered set of $\OCD(C;A)$ for all $\dbinom{m}{n-1}$ unordered $(n-1)$-point subsets $A\subseteq C$. 
The \emph{mirror image} $\mSCD(C)$ is obtained from $\SCD(C)$ by reversing the signs of $(n+1)$ coordinates $s(q;A)$ of $d(q)\in M(C;A)$.
\bs
\end{dfn}

Even if a point of $A$ accidentally coincides with $O(C)$, Definition~\ref{dfn:SCD} still makes sense.
In dimension $n=2$, Definition~\ref{dfn:SCD} needs no permutations because $n-1=1$.
The points of $M(C;\ti A)$ can be lexicographically ordered without affecting future metrics. 
Some of the $\dbinom{m}{n-1}$ $\OCD$s in $\SCD(C)$ can be identical as in Example~\ref{exa:SCD}(b).
If we collapse any $l>1$ identical $\OCD$s into a single $\OCD$ with the \emph{weight} $l/\dbinom{m}{n-1}$, then $\SCD(C)$ can be considered a weighted probability distribution of $\OCD$s.

\begin{exa}[Simplexwise Centered Distribution $\SCD$s for regular clouds]
\label{exa:SCD}
\textbf{(a)}
Let $C_3\subset\R^2$ be the triangular cloud of the points $p_1=(1,0)$, $p_2=(\cos\frac{2\pi}{3},\sin\frac{2\pi}{3})$, $p_3=(\cos\frac{4\pi}{3},\sin\frac{4\pi}{3})$, which form an equilateral triangle with the center of mass $O(C_3)=(0,0)=\vb*{0}$ and side lengths $\sqrt{3}$.
For the base sequence $A_1=\{p_1\}$ and extended $\ti A_1=\{p_1,O(C)\}$, we get 
$$\OCD(C_3;A_1)=\left[|\vb*{p_1}-\vb*{0}|,\left( \begin{array}{cc} 
|\vb*{p_2}-\vb*{p_1}| & |\vb*{p_3}-\vb*{p_1}| \\
|\vb*{p_2}-\vb*{0}|  & |\vb*{p_3}-\vb*{0}|  \\ 
\sign(\vb*{p_2}-\vb*{p_1},\vb*{p_2}) & 
\sign(\vb*{p_3}-\vb*{p_1},\vb*{p_3})
\end{array}\right) \right]=$$
$=\left[1,\left( \begin{array}{cc} 
\sqrt{3} & \sqrt{3} \\
1 & 1  \\ 
-1 & +1
\end{array}\right) \right],$
because
$\det(\vb*{p_2}-\vb*{p_1},\vb*{p_2})=\dettwo{-\frac{3}{2}}{-\frac{1}{2}}{\frac{\sqrt{3}}{2}}{\frac{\sqrt{3}}{2}}<0$ and
$\det(\vb*{p_3}-\vb*{p_1},\vb*{p_2})=\dettwo{-\frac{3}{2}}{-\frac{1}{2}}{-\frac{\sqrt{3}}{2}}{-\frac{\sqrt{3}}{2}}>0$.
Due to the 3-fold rotational symmetry, $\OCD(C_3;A_i)$ is independent of $i=1,2,3$.
Hence, $\SCD(C_3)$ can be considered a distribution of a single $\OCD(C_3;A_1)$ of weight 1.
\medskip

\noindent
\textbf{(b)}
Let $C_4\subset\R^2$ be the square cloud of the points $p_1=(1,0)$, $p_2=(0,1)$, $p_3=(-1,0)$, $p_4=(0,-1)$, see the last image in Fig.~\ref{fig:clouds_oriented}.
The center of mass is $\vb*{0}\in\R^2$ and has a distance $1$ to each $p_i$.
For $A_1=\{p_1\}$, we get
$$M(C_4;\ti A_1)=\left( \begin{array}{ccc} 
|\vb*{p_2}-\vb*{p_1}| & |\vb*{p_3}-\vb*{p_1}| & |\vb*{p_4}-\vb*{p_1}|\\
|\vb*{p_2}-\vb*{0}|  & |\vb*{p_3}-\vb*{0}|   & |\vb*{p_4}-\vb*{0}|  \\ 
\sign(\vb*{p_2}-\vb*{p_1},\vb*{p_2}) & 
\sign(\vb*{p_3}-\vb*{p_1},\vb*{p_3}) & 
\sign(\vb*{p_4}-\vb*{p_1},\vb*{p_4})
\end{array}\right)$$
$=\left( \begin{array}{ccc} 
\sqrt{2} & 2 & \sqrt{2} \\
1 & 1 & 1 \\ 
-1 & 0 & +1
\end{array}\right),$
because
$\det(\vb*{p_2}-\vb*{p_1},\vb*{p_2})=\dettwo{-1}{0}{-1}{1}<0$, then
$\det(\vb*{p_3}-\vb*{p_1},\vb*{p_3})=\dettwo{-2}{-1}{0}{0}=0$, and
$\det(\vb*{p_4}-\vb*{p_1},\vb*{p_4})=\dettwo{-1}{0}{1}{-1}>0$.
Due to the 4-fold rotational symmetry, $\OCD(C_4;A_i)=[1,M(C_4;\ti A_i)]$ is independent of $A_i=\{p_i\}$ for $i=1,2,3,4$.
Hence, $\SCD(C_4)$ can be considered a distribution of a single $\OCD(C_4;A_1)$ of weight 1.
\bs
\end{exa}

Lemma~\ref{lem:OCD_reconstruction} and Theorem~\ref{thm:SCD_complete} follow from Lemma~\ref{lem:ORD_reconstruction} and Theorem~\ref{thm:OSD_complete}, respectively, by adding the center of mass $O(C)$ as the $n$-th point to any base sequence $A\subset C$ of $n-1$ points.

\begin{lem}[reconstruction 
from $\OCD$]
\label{lem:OCD_reconstruction}
Any $n$-dimensional cloud $C\subset\R^n$ of $m$ unordered points with the center of mass $O(C)$ can be reconstructed, uniquely under rigid motion, from $\OCD(C;A)$ in Definition~\ref{dfn:OSD} for any base sequence $A\subseteq C$ of $n-1$ points with $\aff(\ti A)=n-1$.  
\bs
\end{lem}

\begin{thm}[complexity and completeness of $\SCD$ under rigid motion]
\label{thm:SCD_complete}
\textbf{(a)}
For any $n$-dimensional cloud $C$ of $m$ unordered points, the Simplexwise Centered Distribution $\SCD(C)$ is computable in time $O(m^{n})$. 
\medskip

\noindent
\textbf{(b)}
Any $n$-dimensional clouds $C,C'\subset\R^n$ of $m$ unordered points can be exactly matched by rigid motion \emph{if and only if} $\SCD(C)=\SCD(C')$.
\medskip

\noindent
\textbf{(c)}
Any $n$-dimensional clouds $C,C'\subset\R^n$ of $m$ unordered points can be exactly matched by isometry \emph{if and only if} $\SCD(C)$ equals $\SCD(C')$ or its mirror image $\mSCD(C')$.
\bs
\end{thm}

Lemma~\ref{lem:OCD_reconstruction} and Theorem~\ref{thm:SCD_complete} fulfill conditions (a,b) and partially (e), which is the invariant part of Problem~\ref{pro:cloud_invariants_SCD}.
Equality $\SCD(C)=\SCD(C')$ is interpreted as a bijection between the unordered distributions $\SCD(C)\to\SCD(C')$ matching all $\OCD$s, which can be verified by checking if a distance metric between these $\SCD$s equals 0.
The metric part of Problem~\ref{pro:cloud_invariants_SCD} requires a Lipschitz continuous metric on $\SCD$s.
\medskip

Example~\ref{exa:SCD}(b) illustrates the key discontinuity challenge: if $p_3=(-1,0)$ is perturbed, the corresponding sign can discontinuously change to $+1$ or $-1$.
To get a continuous metric on $\SCD$s, we will multiply each sign by a continuous \emph{strength} of a simplex, as formalized below.

\begin{dfn}[a \emph{simplex} $T$ on $n+1$ points in $\R^n$]
\label{dfn:geometric_simplex}
The geometric \emph{simplex} $T$ on any ordered points $p_0,\dots,p_n\in\R^n$ is the subset
$T=\left\{ \sum\limits_{i=0}^n t_i \vb*{p_i} \; \vline \; t_i\in[0,1],\;  \sum\limits_{i=0}^n t_i=1\right\}\subset\R^n$ with ordered \emph{vertices} $p_0,\dots,p_n$.
Define $\sign(T)=\sign(\vb*{p_1}-\vb*{p_0},\dots,\vb*{p_n}-\vb*{p_0})$, see Definition~\ref{dfn:D(A)+M(C;A)}.
\bs
\end{dfn}

The \emph{volume} $\Vol(T)$ of a simplex $T$ or an arbitrary polyhedron is often used as a shape descriptor and detects affine independence in the sense that $T$ is \emph{degenerate} if and only if $\Vol(T)=0$.
However, the volume and all other distance-based descriptors do not distinguish mirror images, which have different signs of orientation.
If a simplex $T$ goes through a degenerate configuration $C$ of vertices with $\aff(C)<n$, then $\sign(T)$ can discontinuously change.
This discontinuity is an obstacle to recognizing nearly mirror-symmetric simplices or more general clouds.
\smallskip

One way to resolve this discontinuity is to consider the signed volume $\sign(T)\Vol(T)$, because $\Vol(T)$ vanishes only on degenerate simplices.
However, the volume $\Vol(T)$ is not Lipschitz continuous, as illustrated below.

\begin{exa}[the area of a triangle is not Lipschitz continuous]
\label{exa:area_not_Lipschitz}
For any large $l>0$ and small $\ep>0$, let $T(l,\ep)\subset\R^2$ be the 2D simplex (triangle) on the vertices $(0,\ep)$ and $(\pm l,0)$.
The signed area $l\ep$ of $T(l,\ep)$ distinguishes mirror images $T(l,\pm\ep)$ but is not Lipschitz continuous under perturbations.
Indeed, as $\ep\to 0$, the triangle $T(l,\ep)$ degenerates to a straight line, while the area drops to 0 too quickly so that $\dfrac{\Vol(T(l,\ep)) - \Vol(T(l,0))}{\ep-0}=\dfrac{l\ep}{\ep}=l$ is not bounded.
Hence, if given points are not restricted to a fixed region, a small change in their positions may lead to a large change for the area of a triangle, and similarly 
in $\R^n$.
\bs
\end{exa}

\index{strength of a simplex}
\begin{dfn}[\emph{strength} $\si(T)$ of a simplex]
\label{dfn:strength_simplex}
Let $T$ be the simplex on any $p_0,\dots,p_n\in\R^n$.
The \emph{half-perimeter} is $p(T)=\dfrac{1}{2}\sum\limits_{i\neq j}|p_i-p_j|$. 
The \emph{strength} is $\si(T)=\dfrac{\Vol^2(T)}{p^{2n-1}(T)}$.
The \emph{signed strength} is $s(T)=\sign(T)\si(T)$. 
\bs
\end{dfn}

\begin{exa}[the strength for $n=1,2$]
\label{exa:strength_n=1,2}
\textbf{(a)}
For $n=1$, the simplex $T$ on ordered points $p_0,p_1\in\R$ is the line segment with the half-perimeter $p(T)=\dfrac{1}{2}|p_1-p_0|$, length $\Vol(T)=|p_1-p_0|=2p(T)$, strength $\si(T)=\dfrac{\Vol^2(T)}{p(T)}=4p(T)=2|p_1-p_0|$, and signed strength $s(T)=2(p_1-p_0)$.
\medskip

\noindent
\textbf{(b)}
For $n=2$, any triangle $T\subset\R^2$ with inter-point distances $a,b,c$ has the half-perimeter
$p=\dfrac{a+b+c}{2}$, area $\Vol(T)=\sqrt{p(p-a)(p-b)(p-c)}$ by Heron's formula, strength $\si(T)=\dfrac{\Vol^2(T)}{p^3}=\dfrac{(p-a)(p-b)(p-c)}{p^2}$.
\bs
\end{exa}

\begin{thm}[strength properties, {\cite[Theorem 2.4]{anosova2026strength}}]
\label{thm:strength_properties}
\textbf{(a)}
The strengths $\si(T)$ and $s(T)$ of any simplex $T\subset\R^n$ are invariant under isometry and rigid motion in $\R^n$, respectively, and can be computed in time $O(n^3)$.
The uniform scaling of $\R^n$ by any factor $c>0$ multiplies $\si(T)$ and $s(T)$ by $c$.
\medskip

\noindent
\textbf{(b)}
There is a constant $\la_n>0$ such that, for any $\ep>0$, if a simplex $Q$ is obtained from a simplex $T\subset\R^n$ by perturbing every vertex of $T$ within its $\ep$-neighborhood, then $|\si(T)-\si(Q)|\leq 2\la_n\ep$ and $|s(T)-s(Q)|\leq 2\la_n\ep$, where $\la_1=2$, $\la_2=\sqrt{3}$, and
$\la_n\leq \dfrac{2^{n+0.5}}{(n!) n^{2n-4}}$ for $n\geq 3$, e.g. $\la_3<0.37$.
\bs
\end{thm} 

\begin{lem}[signs of sequences and simplices]
\label{lem:signs_sequence_simplex}
For any cloud $C\subset\R^n$, a base sequence $A\subset C$ of $n$ ordered points, and any $q\in C\setminus A$, consider the simplex $T$ on the sequence $A\cup\{q\}$, including $q$ as the $(n+1)$-st point.
Then $\sign(T)$ from Definition~\ref{dfn:geometric_simplex} equals the sign $s(q;A)$ from Definition~\ref{dfn:D(A)+M(C;A)}. 
\end{lem}
\begin{proof}
Definition~\ref{dfn:D(A)+M(C;A)} introduced $s(q;A)$ as the sign of the determinant with the columns $\vb*{q}-\vb*{p_1},\dots,\vb*{q}-\vb*{p_n}$.
Then $s(q;A)=\sign(\vb*{q}-\vb*{p_1},\dots,\vb*{q}-\vb*{p_n})=(-1)^n\sign(T')$ for the simplex $T'$ on the ordered $q,p_1,\dots,p_n$, because $\sign(T')=\sign(\vb*{p_1}-\vb*{q},\dots,\vb*{p_n}-\vb*{q})$ by Definition~\ref{dfn:geometric_simplex}.
By moving the vertex $q$ of $T'$ to the last $(n+1)$-st position, we get the simplex $T$ on $A\cup\{q\}$ with $\sign(T)=(-1)^n\sign(T')$, so $s(q;A)=\sign(T)$.
\end{proof}

Theorem~\ref{thm:strength_properties}(b) and 
Lemma~\ref{lem:signs_sequence_simplex} justify replacing the discontinuous sign $s(q;A)$ 
with the  signed strength $s(A\cup\{q\})=s(q;A)\si(A\cup\{q\})$ to define Lipschitz continuous metrics in the next section.

\section{Continuous metrics on complete invariants}
\label{sec:OSD+SCD_metrics}

This section introduces Lipschitz continuous metrics on $\OSD$s and $\SCD$s by using the signed strength of an ordered sequence from Lemma~\ref{lem:signs_sequence_simplex}.
\medskip

For simplicity and similar to the case of a metric space, we assume that a cloud $C\subset\R^n$ is given by a matrix of pairwise Euclidean distances.
If $C$ is given by Euclidean coordinates of points, then any distance requires $O(n)$ computations, and we should add the factor $n$ in all complexities below, which keeps all times polynomial in the number $m$ of points.
\medskip

Recall that Definition~\ref{dfn:Minkowski+Hausdorff+bottleneck}  introduced the bottleneck distance $\BD$ between point clouds of the same size in a metric space.
Definition~\ref{dfn:BMD} interprets $\BD$ as the bottleneck matching distance for a weighted complete bipartite graph, which is needed for more general distributions $\OSD,\SCD$.

\begin{dfn}[Bottleneck Matching Distance $\BMD$]
\label{dfn:BMD}
Let $\Ga$ be a complete bipartite graph with $m$ white vertices and $m$ black vertices so that every white vertex is connected to every black vertex by a single edge $e$ of a weight $w(e)\geq 0$.
A \emph{vertex matching} of the graph $\Ga$ is a collection $E$ of $m$ disjoint edges (with distinct vertices).
The \emph{weight} $W(E)=\max\limits_{e\in E} w(e)$ is the largest weight of an edge in $E$.
The \emph{Bottleneck Matching Distance} $\BMD(\Ga)=\min\limits_E W(E)$ is the minimum weight of a vertex matching. 
\bs
\end{dfn}

The bottleneck matching distance $\BMD(\Ga)$ for a bipartite graph $\Ga$ with $V$ vertices and $E$ edges is computed in time $O(E\sqrt{V})$, see \cite{hopcroft1973n}.

\begin{lem}[metric axioms for the bottleneck matching distance $\BMD$]
\label{lem:BMD_metric}
Let $B=\{B_1,\dots,B_m\}$ and $C=\{C_1,\dots,C_m\}$ be any unordered distributions of $m$ objects with a base metric $d$.
Define the complete bipartite graph $\Ga(B,C)$ with $m$ black vertices (associated with) $B_1,\dots,B_m$ and $m$ white vertices (associated with)  $C_1,\dots,C_m$,
Set the weight of the edge $e$ joining $B_i,C_j$ as $w(e)=d(B_i,C_j)$.
Then $\BMD(\Ga(B,C))$ in Definition~\ref{dfn:BMD} satisfies all metric axioms on unordered distributions of $m$ objects.
\end{lem}
\newcommand{\lemBMDmetric}{
\begin{proof}[\textbf{Proof of Lemma}~\ref{lem:BMD_metric}]
The first axiom says that $\BMD(B,C)=0$ if and only if the weighted distributions $B,C$ are equal in the sense that there is a bijection $g:B\to C$ so that $d(g(B_i),C_i)=0$,  i.e. $g(B_i)=C_i$ due to the first metric axiom for $d$, $i=1,\dots,m$.
Indeed, if the distributions $B,C$ can be matched by a bijection, we get a vertex matching $E$ of $\Ga(B,C)$ whose all edges have weight $0$, then $\BMD(\Ga(B,C))=0$ by Definition~\ref{dfn:BMD}.
Conversely, if $\BMD(\Ga(B,C))=0$, there is a vertex matching $E$ of $\Ga(B,C)$ with all edge weights $0$, which induces a bijection $B\to C$.
\smallskip

The symmetry $\BMD(\Ga(B,C))=\BMD(\Ga(C,B))$ follows from Definition~\ref{dfn:BMD} and the symmetry of $d$.
To prove the triangle inequality 
$$\BMD(\Ga(B,C))+\BMD(\Ga(C,D))\geq \BMD(\Ga(B,D)),$$ let $E_{BC},E_{CD}$ be optimal vertex matchings in the graphs $\Ga(B,C),\Ga(C,D)$, respectively, such that $\BMD(\Ga(B,C))=W(E_{BC})$ and $\BMD(\Ga(C,D))=W(E_{CD})$, see Definition~\ref{dfn:BMD}.
The composition $E_{BD}=E_{BC}\circ E_{CD}$ is a vertex matching in $\Ga(B,D)$, so $W(E_{BD})\geq \BMD(\Ga(B,D))$.
It suffices to prove that $W(E_{BC})+W(E_{CD})\geq W(E_{BD}).$
Let $e_{BD}$ be an edge with a largest weight from $E_{BD}$, so $W(E_{BD})=w(e_{BD})$.
The edge $e_{BD}$ can be considered a union of edges $e_{BC}\in E_{BC}$ and $e_{CD}\in E_{CD}$.
The triangle inequality
$w(e_{BC})+w(e_{CD})\geq w(e_{BD})=W(E_{BD})$ for $d$ implies that $W(E_{BC})+W(E_{CD})\geq W(E_{BD})$ because both terms in the left-hand side are maximized for all edges (not only $e_{BC},e_{CD}$) from $E_{BC},E_{CD}$. 
\end{proof}
}\lemBMDmetric

Definition~\ref{dfn:ORD+OCD_BD+EMD}  will use $\BD$ for clouds of the same size in $\R^{n+1}$ with the base metric $L_\infty$. 
The Earth Mover's Distance $\EMD$ in Definition~\ref{dfn:EMD} will be applied to clouds of any number of equally weighted points in $\R^{n+1}$ with the base metric $L_\infty$. 
Though the $\EMD$ was introduced in \cite{rubner2000earth} for all distributions, which can have different total weights, \cite[appendix]{rubner2000earth} proved that the $\EMD$ is a metric for distributions of the same total weight.
Hence Definition~\ref{dfn:EMD} considers only normalized distributions of total weight $1$.

\begin{dfn}[Earth Mover's Distance on weighted distributions]
\label{dfn:EMD}
Let $B=\{B_1,\dots,B_k\}$ and $C=\{C_1,\dots,C_l\}$ be unordered sets of objects with weights $w(B_i)$, $i=1,\dots,k$, and $w(C_j)$, $j=1,\dots,l$, respectively such that $\sum\limits_{i=1}^k w(B_i)=1=\sum\limits_{j=1}^l w(C_j)$.
Let $d$ be a \emph{base metric} between any objects $B_i$ and $C_j$.
A \emph{flow} from $B$ to $C$ is a $k\times l$ matrix  whose entry $f_{ij}$ represents a \emph{flow} from $B_i$ to $C_j$.
The \emph{Earth Mover's Distance}  
$\EMD(B,C)=\sum\limits_{i=1}^{k} \sum\limits_{j=1}^{l} f_{ij} d(B_i,C_j)$ is minimized for $f_{ij}\in[0,1]$ subject to 
(1) $\sum\limits_{j=1}^{l} f_{ij}\leq w(B_i)$ for any fixed index $i=1,\dots,k$, 
(2) $\sum\limits_{i=1}^{k} f_{ij}\leq w(C_j)$ for any fixed index $j=1,\dots,l$, and
(3) $\sum\limits_{i=1}^{k}\sum\limits_{j=1}^{l} f_{ij}=1$.
\bs
\end{dfn}

The first condition $\sum\limits_{j=1}^{l} f_{ij}\leq w(B_i)$ means that, for each $i=1,\dots,k$, at most the weight $w(B_i)$ `flows' into all objects $C_j$ via $f_{ij}$ for $j=1,\dots,l$. 
The second condition $\sum\limits_{i=1}^{k} f_{ij}\leq w(C_j)$ means that, for each $j=1,\dots,l$, all $f_{ij}$ `flow' from $B_i$ for $i=1,\dots,k$ into $C_j$ up to its weight $w(C_j)$.
The last condition $\sum\limits_{i=1}^{k}\sum\limits_{j=1}^{l} f_{ij}=1$
 forces to `flow' all $B_i$ to all $C_j$.  
The EMD is a partial case of Wasserstein metrics \cite{vaserstein1969markov} in transportation theory \cite{kantorovich1960mathematical}.

\begin{lem}[$\BMD\geq\EMD$]
\label{lem:BMD>=EMD} 
Let $B=\{B_1,\dots,B_m\}$, $C=\{C_1,\dots,C_m\}$ be unordered sets of equally weighted objects with a base metric $d$ and the complete bipartite graph $\Ga(B,C)$ in Lemma~\ref{lem:BMD_metric}.
Fix the weight of $w(e)$ of the edge between any vertices $B_i$ and $C_j$ of $\Ga(B,C)$ as $d(B_i,C_j)$.
Then $\BMD(\Ga(B,C))\geq \EMD(B,C)$, where 
$\EMD(B,C)$ uses any base metric $d'$ such that $d'(B_i,C_j)\leq d(B_i,C_j)$ for all $i,j=1,\dots,m$.
\end{lem}
\begin{proof}
In Definition~\ref{dfn:BMD}, let $E$ be a vertex matching of $\Ga(B,C)$ that minimizes the Bottleneck Matching Distance $\BMD(\Ga(B,C))=\max\limits_{e\in E} w(E)$.
This vertex matching $E$ defines a bijection $g:B\to C$, which is a permutation of $m$ objects.
In the notations of Definition~\ref{dfn:EMD}, $g$ induces the `flow' from $B$ to $C$ with $f_{i,g(i)}=\dfrac{1}{m}$ and $f_{ij}=0$ for all $j\neq g(i)$, $i=1,\dots,m$, with the cost $\sum\limits_{i=1}^{m} \sum\limits_{j=1}^{m} f_{ij} d'(B_i,C_j)=\sum\limits_{i=1}^{m} \dfrac{1}{m}d'(B_i,C_{g(i)})\leq \sum\limits_{i=1}^{m} \dfrac{1}{m}d(B_i,C_{g(i)}) \leq \dfrac{1}{m}\sum\limits_{i=1}^{m} \BMD(\Ga(B,C))=\BMD(\Ga(B,C))$.
Since $\EMD(B,C)$ is the cost minimized over all $f_{ij}\in[0,1]$ subject to the weight-respecting conditions in Definition~\ref{dfn:EMD}, $\BMD(\Ga(B,C))$ is the upper bound of $\EMD(B,C)$.
\end{proof}

In Definition~\ref{dfn:ORD+OCD_BD+EMD}, we use the notations $\BD_\infty$ and $\EMD_\infty$ for metrics on clouds of unordered and equally weighted points with base metric $L_\infty$. 

\begin{dfn}[metrics $\BD_\infty,\EMD_\infty$ on $\ORD,\OCD$]
\label{dfn:ORD+OCD_BD+EMD} 
\textbf{(a)}
Let $C\subset\R^n$ be any cloud of $m$ unordered points with a base sequence $A$ of $n$ ordered points.
In the notations of Definition~\ref{dfn:OSD}, for each point $q\in C\setminus A$, replace the sign $s(q;A)$ in the last coordinate of $d(q)\in\R^{n+1}$ with $\dfrac{s(A\cup\{q\})}{\la_n}$, where $s(A\cup\{q\})=s(q;A)\si(A\cup\{q\})$ is called the \emph{signed strength}.
\smallskip

Denote the resulting cloud of $m-n$ unordered points $\hat d(q)\in\R^{n+1}$ by $\hat M(C;A)$.
Any permutation $\xi\in S_n$ of $1,\dots,n$ acts on $D(A)$ and $\hat M(C;A)$ as in Definition~\ref{dfn:OSD}.
For another cloud $C'\subset\R^n$ of unordered points with a base sequence $A'$ of $n$ ordered points, similarly define the cloud $\hat M(C';A')$.
Let the cloud $C'$ have the same number of points as $C$.
Then set 
$$\BD_\infty[\xi]=\max\Big\{\; L_{\infty}\big(\xi(D(A)),D(A')\big),\; \BD_\infty\big(\xi(\hat M(C;A)),\hat M(C';A')\big) \; \Big\}.$$
In general, even if the sizes of the clouds $C,C'\subset\R^n$ differ, set
$$\EMD_\infty[\xi]=\max\big\{ L_{\infty}\big(\xi(D(A)),D(A')\big), \EMD_\infty\big(\xi(\hat M(C;A)),\hat M(C';A')\big) \big\}.$$
Define the following metrics between Oriented Relative Distributions:
$$\begin{array}{l}
\BD_\infty\big(\ORD(C;A),\ORD(C';A')\big)=\min\limits_{\xi\in S_n} \BD_\infty[\xi], \text{ and}\\
\EMD_\infty\big(\ORD(C;A),\ORD(C';A')\big)=\min\limits_{\xi\in S_n} \EMD_\infty[\xi].
\end{array}$$

\noindent
\textbf{(b)}
The metrics $\BD_\infty,\EMD_\infty$ are defined on $\OCD$s similarly to part (a) after replacing $A\subset C$ with the extended base sequence $\ti A=A\cup\{O(C)\}$, which includes the center of mass $O(C)$ of a given cloud $C\subset\R^n$.
\bs 
\end{dfn}

In Definition~\ref{dfn:ORD+OCD_BD+EMD} , the subscript $\infty$ refers to the base Chebyshev metric $L_\infty$, which can be replaced with many other base metrics on distance matrices and clouds in $\R^{n+1}$.
For the $(n+1)$-st coordinate of points in clouds $M(C;A)\subset\R^{n+1}$, the coefficient $\dfrac{1}{\la_n}$ in front of the signed strength $s(A\cup\{q\})$ normalises the Lipschitz constant $2\la_n$ from Theorem~\ref{thm:strength_properties}(b) to $2$ in line with changes of inter-point distances by at most $2\ep$ when given points are perturbed within their $\ep$-neighbourhoods.

\begin{lem}[properties of metrics $\BD_\infty$ and $\EMD_\infty$ on $\ORD$s]
\label{lem:ORD+OCD_BD+EMD}
\textbf{(a)}
For any $n$-dimensional clouds $C,C'\subset\R^n$ of $m$ unordered points with base sequences $A\subset C$ and $A'\subset C'$ of $n$ ordered points, the bottleneck distance $\BD_\infty\big(\ORD(C;A),\ORD(C';A')\big)$ in Definition~\ref{dfn:ORD+OCD_BD+EMD}(a) satisfies all metric axioms and can be computed in time $O\big(n!(n^2+m^{1.5}\log^{n+1} m)\big)$.
Also,
$\BD_\infty\big(\ORD(C;A),\ORD(C';A')\big)\geq 
\EMD_\infty\big(\ORD(C;A),\ORD(C';A')\big)$.
\medskip

\noindent
\textbf{(b)}
For any $n$-dimensional clouds $C,C'\subset\R^n$ of up to $m$ unordered points with base sequences $A\subset C$ and $A'\subset C'$ of $n$ ordered points, $\EMD_\infty\big(\ORD(C;A),\ORD(C';A')\big)$ in Definition~\ref{dfn:ORD+OCD_BD+EMD}(a) satisfies all metric axioms and can be computed in time $O\big(n!(n^2+m^3\log m)\big)$.
\end{lem}
\begin{proof}
The first metric axiom says that $\ORD(C;A),\ORD(C';A')$ are equivalent by Definition~\ref{dfn:OSD} if and only if $\BD_{\infty}\big(\ORD(C;A),\ORD(C';A')\big)=0$, i.e. $\BD_\infty[\xi]=0$ for some permutation $\xi\in S_n$.
Due to the first metric axiom for $L_\infty$ and $\BD$, the equality $\BD_\infty[\xi]=0$ is equivalent to $\xi(D(A))=D(A')$ and $\xi(M(C;A))=M(C';A')$ as clouds of unordered points in $\R^{n+1}$.
The last two conclusions mean that the Oriented Relative Distributions $\ORD(C;A),\ORD(C';A')$ are equivalent by Definition~\ref{dfn:OSD}.
\smallskip

The symmetry axiom for $\BD_\infty$ follows because any permutation $\xi\in S_n$ is invertible.
To prove the triangle inequality 
$$\begin{array}{l}
\BD_{\infty}\big(\ORD(C;A),\ORD(C';A')\big)+
\BD_{\infty}\big(\ORD(C';A'),\ORD(C'';A'')\big)\geq \\
\BD_{\infty}\big(\ORD(C;A),\ORD(C'';A'')\big),
\end{array}$$ let $\xi,\xi'\in S_n$ be optimal permutations for the $\BD_\infty$ values in the left-hand side above. 
The triangle inequality for $L_\infty$ says that 
$$\begin{array}{l}
L_{\infty}\big(\xi(D(A)),D(A')\big)+
L_{\infty}\big(\xi'(D(A'')),D(A')\big)\geq \\
L_{\infty}\big(\xi(D(A)),\xi'(D(A''))\big)=
L_{\infty}\big((\xi')^{-1}\xi(D(A)),D(A'')\big),
\end{array}$$ similarly for the bottleneck distance $\BD$.
Taking the maximum of $L_\infty,W_\infty$ preserves the triangle inequality, see other metric transforms in \cite[section 4.1]{deza2009encyclopedia}.
Then $\BD_{\infty}\big(\ORD(C;A),\ORD(C'' ;A'')\big)=\min\limits_{\xi\in S_n}d(\xi)$ cannot be larger than $\BD_{\infty}[(\xi')^{-1}\xi]$ for the composition of the permutations above, so the triangle inequality holds for $\BD_\infty$. 
For a fixed $\xi\in S_{n}$, the distance $L_\infty(\xi(D(A)),D(A'))$ requires time $O(n^2)$.
By Theorem~\ref{thm:strength_properties}(a) the signed strengths in the $(n+1)$-st coordinates of all points of $\hat M(C;A),\hat M(C';A')$ need time $O(mn^3)$, which is dominated by the final time. 
The bottleneck distance $\BD_{\infty}$ between the clouds $\xi(\hat M(C;A))$ and $\hat M(C';A')$ of $m-n$ unordered points in $\R^{n+1}$ needs time $O(m^{1.5}\log^{n+1} m)$ by
\cite[Theorem~6.5]{efrat2001geometry}.
The minimization over all $\xi\in S_n$ gives the factor $n!$ in the final time.
\smallskip
 
To prove that $\BD_\infty\geq\EMD_\infty$, take the permutation $\xi\in S_n$ that minimizes $\BD_\infty[\xi]$, i.e. $\BD_\infty\big(\ORD(C;A),\ORD(C';A')\big)=\BD_\infty[\xi]$.
Since both metrics $\BD_\infty$ and $\EMD_\infty$ between clouds use the base distance $L_\infty$ between points in $\R^{n+1}$, Lemma~\ref{lem:BMD>=EMD} implies that
$$\BD_\infty\big(\xi(\hat M(C;A)),\hat M(C';A')\big)\geq 
\EMD_\infty\big(\xi(\hat M(C;A)),\hat M(C';A')\big).$$
Taking the maximum of both sides with $L_{\infty}\big(\xi(D(A)),D(A')\big)$ gives
$$\BD_\infty\big(\ORD(C;A),\ORD(C';A')\big)=\BD_\infty[\xi]\geq\EMD_\infty[\xi].$$
Taking the minimum of the right hand side over $\xi\in S_n$ gives
$$\BD_\infty\big(\ORD(C;A),\ORD(C';A')\big)\geq
\EMD_\infty\big(\ORD(C;A),\ORD(C';A')\big).$$

\noindent
\textbf{(b)}
The metric axioms for $\EMD_\infty$ are proved as in  part (a) by replacing $\BD$ with $\EMD$ whose metric axioms were checked in \cite[appendix]{rubner2000earth}.
Consider the clouds $\hat M(C;A)$ and $\hat M(C';A')$ as unordered distributions of equally weighted $m-n$ points in $\R^{n+1}$.
The $\EMD$ between weighted distributions of a maximum size $m$ can be computed in time $O(m^3\log m)$, see 
\cite{fredman1987fibonacci,goldberg1987solving}.
Taking into account the time $O(n^2)$ for $L_\infty(\xi(D(A)),D(A'))$ and minimizing over $\xi\in S_n$ gives $O\big(n!(n^2+m^3\log m)\big)$.
\end{proof}

Definition~\ref{dfn:OSD+SCD_BD+EMD} adapts $\BD_\infty,\EMD_\infty$ in Definition~\ref{dfn:ORD+OCD_BD+EMD}  to 
invariants. 
 
\begin{dfn}[metrics $\BD_\infty$ and $\EMD_\infty$ on invariants $\OSD$ and $\SCD$]
\label{dfn:OSD+SCD_BD+EMD} 
\textbf{(a)}
Let $C,C'\subset\R^n$ be any $n$-dimensional clouds of $m$ unordered points.
Let $\Ga(C,C')$ denote the complete bipartite graph with $\dbinom{m}{n}$ black vertices and $\dbinom{m}{n}$ white vertices that are associated with $\ORD(C;A)$ and $\ORD(C';A')$ for all base sequences $A\subset C$ and $A'\subset C'$ of $n$ ordered points.
For each edge of $\Ga(C,C')$ between vertices associated $\ORD(C;A)$ and $\ORD(C';A')$, define its weight as $\BD_\infty\big(\ORD(C;A),\ORD(C';A')\big)$ in Definition~\ref{dfn:ORD+OCD_BD+EMD} (a).
The \emph{bottleneck distance} $\BD_\infty\big(\OSD(C),\OSD(C')\big)$ is the bottleneck matching distance $\BMD(\Ga(C,C'))$.
In general, even if the sizes of the clouds $C,C'$ differ, we can consider $\OSD(C)$ and $\OSD(C')$ as weighted distributions whose elements $\ORD(C;A)$ and $\ORD(C';A')$ have the base metric $\EMD_\infty\big(\ORD(C;A),\ORD(C';A')\big)$.
Define the \emph{Earth Mover's Distance} $\EMD_\infty\big(\OSD(C),\OSD(C')\big)$ as in Definition~\ref{dfn:EMD}.
\medskip

\noindent
\textbf{(b)}
The metrics $\BD_\infty,\EMD_\infty$ are defined on $\SCD$s similarly to part (a) after replacing $A\subset C$ with the extended base sequence $\ti A=A\cup\{O(C)\}$, which includes the center of mass $O(C)$ of a given cloud $C\subset\R^n$.
\bs
\end{dfn}

\begin{thm}[properties of $\BD_\infty$ and $\EMD_\infty$ on $\OSD$s]
\label{thm:OSD_metrics}
\textbf{(a)}
For any $n$-dimensional clouds $C,C'$ of $m$ unordered points, 
 $\BD_\infty\big(\OSD(C),\OSD(C')\big)$ in Definition~\ref{dfn:OSD+SCD_BD+EMD}(a) satisfies all metric axioms and can be computed in time $O\big(n!(n^2+m^{1.5}\log^{n+1} m)k^2+k^{2.5}\big)$ for $k=\dbinom{m}{n}=O(m^n)$.
Also, 
$\BD_\infty(\OSD(C),\OSD(C'))\geq \EMD_\infty(\OSD(C),\OSD(C'))$.
\medskip

\noindent
\textbf{(b)}
Let $C,C'\subset\R^n$ be any $n$-dimensional clouds, which may consist of up to $m$ unordered points.
Let $\OSD(C),\OSD(C')$ have a maximum size  $l\leq\dbinom{m}{n}$ after collapsing identical  $\ORD$s.
Then $\EMD_\infty\big(\OSD(C),\OSD(C')\big)$ in Definition~\ref{dfn:OSD+SCD_BD+EMD}(a) satisfies all metric axioms and can be computed in time $O\big(n!(n^2 +m^{1.5}\log^{n+1} m)l^2+ l^3\log l\big)$.
\end{thm}
\begin{proof}
\textbf{(a)}
The base metrics $\BD_\infty,\EMD_\infty$ on $\ORD$s in Definition~\ref{dfn:ORD+OCD_BD+EMD}(a) satisfy the axioms by Lemma~\ref{lem:ORD+OCD_BD+EMD}(a).
Then $\BD_\infty,\EMD_\infty$ on $\OSD$s in Definition~\ref{dfn:OSD+SCD_BD+EMD}(a) satisfy the axioms by Lemma~\ref{lem:BMD_metric}, \cite{rubner2000earth}, respectively.
\smallskip

In part (a), both clouds $C,C'\subset\R^n$ consist of $m$ unordered points. 
Then both $\OSD(C),\OSD(C')$ consist of $k=\dbinom{m}{n}$ distributions $\ORD$.
The graph $\Ga(C,C')$ in Definition~\ref{dfn:BMD} has $V=2k$ vertices and $E=k^2$ edges whose all weights are computed in time $O\big(n!(n^2+m^{1.5}\log^{n+1} m)k^2\big)$.
Then $\BD_\infty\big(\OSD(C),\OSD(C')\big)=\BMD\big(\Ga(C,C')\big)$ is computed by \cite{hopcroft1973n} in extra time $O(E\sqrt{V})=O(k^{2.5})$.
Lemma~\ref{lem:BMD>=EMD} and the inequality $\BD_\infty\geq\EMD_\infty$ on $\ORD$s in Lemma~\ref{lem:ORD+OCD_BD+EMD}(a) imply that $\BD_\infty\geq\EMD_\infty$ holds for $\OSD$s.
\medskip

\noindent
\textbf{(b)}
$\OSD(C),\OSD(C')$ consist of $l\leq \dbinom{m}{n}$ distributions $\ORD$.
Each $\EMD_\infty$ between all $\ORD$s in the invariants $\OSD(C)$ and $\OSD(C')$ is computable in time $O\big(n!(n^2 +m^{1.5}\log^{n+1} m)l\big)$.
The total time of all $E\leq l^2$ distances is 
$O\big(n!(n^2 +m^{1.5}\log^{n+1} m)l^2\big)$.
The final extra time $O(l^3\log l)$ is for computing $\EMD$ between distributions whose size is at most $l$. 
\end{proof}

Theorem~\ref{thm:SCD_metrics} is proved as Theorem~\ref{thm:OSD_metrics} by replacing $n$ with $n-1$, because any base sequence $A\subset C$ of $n-1$ ordered points is extended by the center mass $O(C)$ to the sequence $\ti A=A\cup\{O(C)\}$ of $n$ points.

\begin{thm}[properties of $\BD_\infty$ and $\EMD_\infty$ on invariants $\SCD$s]
\label{thm:SCD_metrics}
\textbf{(a)}
For any $n$-dimensional clouds $C,C'\subset\R^n$ of $m$ unordered points, the bottleneck distance $\BD_\infty\big(\SCD(C),\SCD(C')\big)$ in Definition~\ref{dfn:OSD+SCD_BD+EMD}(b) satisfies all metric axioms and,  for $k=\dbinom{m}{n-1}=O(m^{n-1})$, can be computed in time $O\big((n-1)!(n^2+m^{1.5}\log^{n} m)k^2+k^{2.5}\big)$.
Also, the inequality
$\BD_\infty(\SCD(C),\SCD(C'))\geq \EMD_\infty(\SCD(C),\SCD(C'))$ holds.
\medskip

\noindent
\textbf{(b)}
Let $C,C'\subset\R^n$ be any $n$-dimensional clouds, which can have up to $m$ unordered points.
Let $\SCD(C),\SCD(C')$ have a maximum size  $l\leq\dbinom{m}{n-1}$ after collapsing identical  $\ORD$s.
Then $\EMD_\infty\big(\SCD(C),\SCD(C')\big)$ in Definition~\ref{dfn:OSD+SCD_BD+EMD}(b) satisfies all metric axioms and can be computed in time $O\big((n-1)!(n^2 +m^{1.5}\log^{n} m)l^2+ l^3\log l\big)$.
\bs
\end{thm}

A rough upper bound is $l\leq k=\dbinom{m}{n-1}=\dfrac{m(m-1)\dots(m-n+2)}{n!}=O(m^{n-1})$.
In dimension $n=2$, Theorem~\ref{thm:SCD_metrics} gives time $O(m^4\log^2 m)$ for $\BD_\infty$ and $\EMD_\infty$.
This upper bound is faster than the past time $O(m^5\log m)$ for comparing $m$-point clouds by the Hausdorff distance minimized over all isometries in $\R^2$ \cite{chew1997geometric}.
For $n\geq 3$, $\BD_\infty$ and $\EMD_\infty$ are computable in times $O(m^{2.5(n-1)})$ and $O(nm^{3(n-1)}\log m)$, respectively.

\index{strength of a simplex}
\index{Lipschitz continuity}
\begin{thm}[Lipschitz continuity of $\OSD,\SCD$]
\label{thm:OSD+SCD_continuous}
For any cloud $C\subset\R^n$ of $m$ unordered points, perturbing any point within its $\ep$-neighborhood changes $\OSD(C)$ and $\SCD(C)$ by at most $2\ep$ in $\BD_\infty$ and $\EMD_\infty$.
\end{thm}
\begin{proof}
Arbitrarily order points of $C$ as $p_1,\dots,p_m$.
Then $C'$ consist of their $\ep$-perturbations $p'_1,\dots,p'_m$ such that $|p_i-p'_i|\leq\ep$ for $i=1,\dots,m$.
\smallskip

The triangle inequality implies that any inter-point distance changes by at most $2\ep$, i.e.
$\big| |p_i-p_j| - |p'_i-p'_j| \big|\leq |p_i-p'_i|+|p_j-p'_j|\leq 2\ep$ for any $i,j=1,\dots,m$.
Without loss of generality, one can assume that a base sequence $A\subset C$ consists of the first $n$ points $p_1,\dots,p_n$.
\smallskip

Let $\xi\in S_n$ be the trivial permutation keeping all $1,\dots,n$ fixed. 
Let $A'\subset C'$ denote the base sequence of perturbed points $p'_1,\dots,p'_n$.
In the notations of Definition~\ref{dfn:ORD+OCD_BD+EMD}, any element (coordinate) of $D(A)=\xi(D(A))$ and $\hat M(C;A)=\xi(\hat M(C;A))$ differs from the corresponding element of $D(A')$ and $\hat M(C';A')$ by at most $2\ep$.
This upper bound of $2\ep$ for normalized strengths, which are the last $(n+1)$-st coordinates of points in $\hat M(C;A)$ and $\hat M(C';A')$, follows from Theorem~\ref{thm:strength_properties}(b) and Lemma~\ref{lem:signs_sequence_simplex}.
Then all three distances $L_\infty(\xi(D(A)),D(A'))$, $\BD_\infty(\xi(\hat M(C;A)),\hat M(C';A'))$, and $\EMD_\infty(\xi(\hat M(C;A)),\hat M(C';A'))$ have the upper bound $2\ep$.
Minimizing over all permutations $\xi\in S_n$ preserves the upper bound $2\ep$, i.e.
$$\begin{array}{l}
\BD_\infty\big(\ORD(C;A),\ORD(C';A')\big)=\min\limits_{\xi\in S_n} \BD_\infty[\xi]\leq 2\ep, \\
\EMD_\infty\big(\ORD(C;A),\ORD(C';A')\big)=\min\limits_{\xi\in S_n} \EMD_\infty[\xi]\leq 2\ep.
\end{array}$$
In the notations Definition~\ref{dfn:OSD+SCD_BD+EMD},
the trivial permutation $\xi\in S_n$ induces the matching $E$ for $k=\dbinom{m}{n}$ pairs of vertices of the bipartite graph $\Ga(C,C')$, which are associated with all base sequences $A\subset C$ and their perturbations $A'\subset C'$.
Since all weights of edges $e\in E$ have weights at most $2\ep$, we get $\BD_\infty(\OSD(C),\OSD(C'))=\BMD(\Ga(C,C'))\leq 2\ep$.
Similarly, the trivial flows $\OSD(C)\to\OSD(C')$ with base distances at most $2\ep$ imply that normalized distributions have $\EMD_\infty(\OSD(C),\OSD(C'))\leq 2\ep$.
\smallskip

The arguments above work similarly work for the $\SCD$, because 
the center of mass changes under $\ep$-perturbations by at most $\ep$ as follows:
$|O(C)-O(C')|=\left|\dfrac{1}{m}\sum\limits_{i=1}^m p_i -  \dfrac{1}{m}\sum\limits_{i=1}^m p'_i \right|\leq \dfrac{1}{m}\sum\limits_{i=1}^m |p_i-p'_i|\leq \dfrac{1}{m}\sum\limits_{i=1}^m \ep =\ep$.
\end{proof}

Lemma~\ref{lem:OCD_reconstruction} and
Theorems~\ref{thm:SCD_complete}, \ref{thm:SCD_metrics}, and \ref{thm:OSD+SCD_continuous} imply that the Simplexwise Centered Distribution fully solves Problem~\ref{pro:cloud_invariants_SCD}.

\section{Faster moment invariants and distances }
\label{sec:OSD+SCD_moments}

This section simplifies the complete $\SCD$ to faster invariants based on averages and other statistical moments of unordered distributions.
\smallskip

Let $S$ consist of $s_1,\dots,s_m\in\R$ with weights $w_1,\dots,w_m$, respectively, such that $\sum\limits_{i=1}^m w_i=1$.
For any integer $t\geq 1$, the $t$-th \emph{moment} \cite[section~2.7]{keeping1995introduction} is 
$\mu_t\{S\}=\sqrt[t]{m^{1-l}\sum\limits_{i=1}^m w_i s_i^t}$, so $\mu_1\{S\}=\sum\limits_{i=1}^m w_i s_i$ is the average.

\begin{dfn}[Average Simplexwise Moments $\ASM(C;t)$]
\label{dfn:ASV+ASD+ASM}
For any  cloud $C\subset\R^n$ of $m$ unordered points, let $A\subset C$ be a base sequence of $n$ points $p_1,\dots,p_n$. 
Recall that $\SPD(A)$ from Definition~\ref{dfn:SPD} consists of $\dfrac{n(n-1)}{2}$ inter-point distances of $A$, written in increasing order and hence independent of an order of points in $A$. 
The vector $\ar{M}(C;A)$ is formed by $m-n$ averages $\dfrac{1}{n}\sum\limits_{i=1}^{n}|\vb*{q}-\vb*{p_i}|$ for all points $q\in C\setminus A$, written in increasing order of values.
The vector $\ar{S}(C;A)$ is formed by $m-n$ normalized strengths $\dfrac{\si(A\cup\{q\})}{\la_n}$, written in increasing order, for all points $q\in C\setminus A$.
\medskip

\noindent
\textbf{(a)}
In the notations above, the \emph{Average Simplexwise Vector} is the tuple
$$\ASV(C;A)=[\SPD(A);\ar{M}(C;A);\ar{S}(C;A)] \text{ of }
\dfrac{n(n-1)}{2}+2(m-n) \text{ values}.$$
For any other cloud $C'\subset\R^n$ 
with a base sequence $A'$, define the metric $\EMD_\infty\big(\ASV(C;A),\ASV(C';A')\big)$ as the maximum of 
$L_\infty\big(\SPD(A),\SPD(A')\big)$,
$\EMD_\infty\big(\ar{M}(C;A), \ar{M}(C';A')\big)$, and 
$\EMD_\infty\big(\ar{S}(C;A), \ar{S}(C';A')\big).$
\medskip

\noindent
\textbf{(b)}
The \emph{Average Simplexwise Distribution} $\ASD(C)$ is the unordered set of equally weighted vectors $\ASV(C;A)$ for all $\dbinom{m}{n}$ unordered subsets $A\subseteq C$ of $n$ points.
Use the base metric $\EMD_\infty$ on $\ASV$s in part (a) to define
the Earth Mover's Distance $\EMD_\infty\big(\ASD(C),\ASD(C')\big)$ as in  Definition~\ref{dfn:EMD}.
\medskip

\noindent
\textbf{(c)}
For $t\geq 1$, the \emph{Average Simplexwise Moment} $\ASM(C;t)=\mu_t\{\ASD(C)\}$ is the $t$-th moment of the distribution $\ASD(C)$ of equally weighted vectors $\ASV(C;A)$ over all $\dbinom{m}{n}$ unordered subsets $A\subseteq C$ of $n$ points, i.e.
$$\ASM(C;t)=\big[\mu_t\{\SPD(A)\};\mu_t\{\ar{M}(C;A)\};\mu_t\{\ar{S}(C;A)\}\big].$$
For any clouds $C,C'\subset\R^n$, define $\EMD_\infty\big(\ASM(C;t),\ASM(C';t)\big)$ as the maximum of the three metrics $\EMD_\infty$ between constituent subvectors.
\bs
\end{dfn}

Since $\si(A\cup\{q\})$ is independent of an order of points in $A$, the vector $\ASV(C;A)$ is well-defined for any subset $A\subset C$ of $n$ unordered points.
\smallskip

The metrics $\EMD_\infty$ in Definition~\ref{dfn:ASV+ASD+ASM} make sense for any clouds $C,C'\subset\R^n$, which might have different sizes.
If underlying distributions are 1-dimensional as for Average Simplexwise Moments $\ASM(C;t)$, then $\EMD_\infty$ between them can be analytically computed by Proposition~\ref{prop:EMD_sets1D} below.
 
\begin{prop}[$\EMD$ on weighted distributions in $\R$] 
\label{prop:EMD_sets1D}
Let $S\subset\R$ be a weighted distribution of $m$ unordered numbers with a total weight 1.
Let us write all values of $S$ as a sequence $\ar{S}$ of numbers $s_1\leq\dots\leq s_m$ with weights $w_1,\dots,w_m$, respectively.
Define the \emph{Cumulative Distribution Function} 
$\CDF[S](x)=\left\{\begin{array}{l} 
0 \text{ for } x<s_1; \\
\sum\limits_{i=1}^k w_i \text{ for } s_k\leq x<s_{k+1}, k=1,\dots,m-1; \\
1 \text{ for } x\geq s_m.
\end{array}  \right. $
\smallskip

Any other distribution $S'$ of a total weight 1 satisfies
$\EMD_\infty(S,S')=\int\limits_{\R} \big|\CDF[S](x)-\CDF[S'](x)\big|dx$.
If $S,S'$ consist of the same number of equally weighted values, then $\EMD_\infty(S,S')=\dfrac{1}{m}L_1(\ar{S},\ar{S}')$. 
\end{prop}
\begin{proof}
The equality $\EMD_\infty(S,S')=\int\limits_{\R} \big|\CDF[S](x)-\CDF[S'](x)\big|dx$ is \cite[Theorem~3]{cohen1997earth}.
Let the distributions $S,S'$ consist of the same number $m$ of equally weighted (not necessarily distinct) values $s_1\leq\dots\leq s_m$ and $s'_1\leq\dots\leq s'_m$.
Then the area $\int\limits_{\R} \big|\CDF[S](x)-\CDF[S'](x)\big|dx$ between the graphs $y=\CDF[S](x)$ and $y=\CDF[S'](x)$ can be split by the lines $y=\dfrac{k}{m}$ for $k=1,\dots,m-1$ into $m$ horizontal slices.
For each $k=1,\dots,m$, in the slice right below the line $y=\dfrac{k}{m}$, the difference of areas is $\dfrac{|s_k-s'_k|}{m}$.
Hence the total difference of areas is $\sum\limits_{k=1}^m\dfrac{|s_k-s'_k|}{m}=\dfrac{1}{m}L_1(\ar{S},\ar{S}')$.
\end{proof}

\begin{exa}[$\ASD$ and $\ASM$ for the clouds $T,K\subset\R^2$ in Fig.~\ref{fig:clouds_oriented}]
\label{exa:ASV+ASD+ASM}
For the 4-point clouds $T,K$ in Fig.~\ref{fig:clouds_oriented}, a base sequence $A$ consists of two (unordered) points.
To compute $\dfrac{\si(A\cup\{q\})}{\la_2}$, we use the Lipschitz constant $\la_2=\sqrt{3}$ and the strength of a triangle on $A\cup\{q\}$ in Example~\ref{exa:strength_n=1,2}(b).
\smallskip

The trapezoid $T$ in Fig.~\ref{fig:cloud_T_strengths} has triangles $\triangle p_1 p_2 p_3$ and $\triangle p_1 p_2 p_4$ with sides $\sqrt{2},2,\sqrt{10}$, $\Vol=1$, half-perimeter $p=\dfrac{\sqrt{2}+2+\sqrt{10}}{2}$, and normalized strength $\si_1
=\dfrac{\Vol^2}{\la_2 p^3}=\dfrac{1}{\sqrt{3}}\left(\dfrac{2}{\sqrt{2}+2+\sqrt{10}}\right)^3\approx 0.016$.
\smallskip

The kite $K$ in Fig.~\ref{fig:cloud_K_strengths} has $\triangle p_1 p_2 p_3$ with sides $\sqrt{2},\sqrt{2},2$, $\Vol=1$, half-perimeter $p=1+\sqrt{2}$, normalized strength $\si_3=\dfrac{1}{\sqrt{3}(1+\sqrt{2})^3}\approx 0.041$.

\begin{figure}[h!]
\centering
\includegraphics[width=\linewidth]{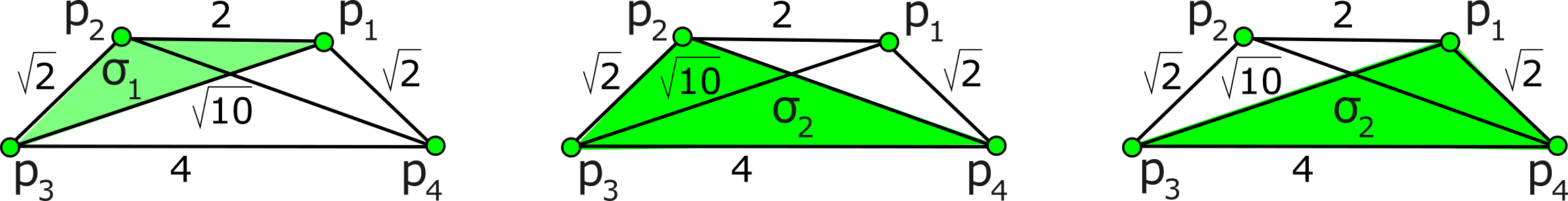}
\caption{The trapezoid $T\subset\R^2$ on points $p_1=(1,1)$, $p_2=(-1,1)$, $p_3=(-2,0)$, $p_4=(2,0)$ has triangles with normalized strengths $\si_1=\dfrac{1}{\sqrt{3}}\left(\dfrac{2}{\sqrt{2}+2+\sqrt{10}}\right)^3\approx 0.016$ and $\si_2=\dfrac{4}{\sqrt{3}}\left(\dfrac{2}{\sqrt{2}+\sqrt{10}+4}\right)^3\approx 0.029$, see Example~\ref{exa:ASV+ASD+ASM}.
}
\label{fig:cloud_T_strengths}
\end{figure}

\begin{figure}[h!]
\centering
\includegraphics[width=\linewidth]{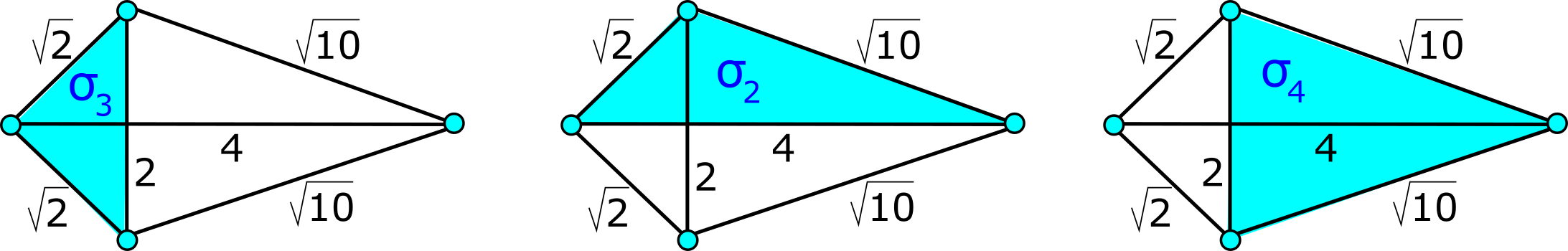}
\caption{The kite cloud $K\subset\R^2$ on points $p_1=(0,1)$, $p_2=(-1,0)$, $p_3=(0,-1)$, $p_4=(3,0)$ has triangles with normalized strengths $\si_3=\dfrac{1}{\sqrt{3}(1+\sqrt{2})^3}\approx 0.041$, $\si_2=\dfrac{2^2}{\sqrt{3}}\left(\dfrac{2}{\sqrt{2}+\sqrt{10}+4}\right)^3\approx 0.029$, $\si_4=\dfrac{3\sqrt{3}}{(1+\sqrt{10})^3}\approx 0.072$.}
\label{fig:cloud_K_strengths}
\end{figure}

Both $T,K$ have triangles with sides $\sqrt{2},\sqrt{10},4$, $\Vol=2$, half-perimeter $p=\dfrac{\sqrt{2}+\sqrt{10}+4}{2}$, and normalized strength
$\si_2=\dfrac{2^2}{\sqrt{3}}\left(\dfrac{2}{\sqrt{2}+\sqrt{10}+4}\right)^3\approx 0.029$.
Finally, $K$ has $\triangle p_1 p_3 p_4$ with sides $2,\sqrt{10},\sqrt{10}$, $\Vol=3$, $p=1+\sqrt{10}$ and normalized strength $\si_4=\dfrac{3^2}{\sqrt{3}}\left(\dfrac{1}{1+\sqrt{10}}\right)^3=\dfrac{3\sqrt{3}}{(1+\sqrt{10})^3}\approx 0.072$.
\medskip

Table~\ref{tab:ASD_TK} summarizes Average Simplexwise Distributions $\ASD(T),\ASD(K)$, each consisting of 6 Average Simplexwise Vectors $\ASV(T;A)$ and $\ASV(K;A)$, respectively, of $1+2+2$ coordinates.
The first coordinate is the distance between points of an unordered subset $A$.
The two further coordinates are the average distances $\bar d(q)$ from Definition~\ref{dfn:ASV+ASD+ASM}  in increasing order.
The two last coordinates are normalized strengths, also in increasing order.
\smallskip

In the trapezoid $T$, the subset $A=\{p_1,p_4\}$ has $\SPD(A)$ consisting of one distance $|p_1-p_4|=\sqrt{2}$.
The average distance to $p_1$ from $p_2,p_3$ is $\dfrac{2+\sqrt{10}}{2}$.
The average distance to $p_4$ from $p_2,p_3$ is $\dfrac{\sqrt{10}+4}{2}$.
The strengths of $\triangle p_1 p_2 p_4$ and $\triangle p_1 p_3 p_4$ are $\si_0<\si_2$, respectively.
The resulting Average Simplexwise Vector $\ASV(T;\{p_1,p_4\})\in\R^5$ is in the first cell of Table~\ref{tab:ASD_TK}.
Another subset $A=\{p_2,p_3\}$ gives the same vector  $\ASV(T;\{p_1,p_4\})$.
This fact is marked by the `double' factor $\times 2$ in front of vectors in Table~\ref{tab:ASD_TK}.
\smallskip

\begin{table*}[h!]
\begin{tabular}{l|l}
$\cbox{yellow}{T}$: Average Simplexwise Vectors  & 
$\cbox{yellow}{K}$: Average Simplexwise Vectors  \\
\hline   

$2\times[\sqrt{2}; \cbox{yellow}{\frac{2+\sqrt{10}}{2}}, \frac{\sqrt{10}+4}{2}; \cbox{yellow}{\si_1},\cbox{yellow}{\si_2} ]\approx$ &
$2\times[\sqrt{2}; \cbox{yellow}{\frac{2+\sqrt{2}}{2}}, \frac{\sqrt{10}+4}{2}; \cbox{yellow}{\si_2}, \cbox{yellow}{\si_3}]\approx$
	\smallskip 	\\

$[1.414; \cbox{yellow}{2.581}, 3.581;  \cbox{yellow}{0.016},  \cbox{yellow}{0.029} ]$ &     
$[1.414; \cbox{yellow}{1.707}, 3.581; \cbox{yellow}{0.029},  \cbox{yellow}{0.041} ]$  
\medskip \\ \hline

$[2;\cbox{yellow}{\frac{\sqrt{2}+\sqrt{10}}{2}}, \; \cbox{yellow}{\frac{\sqrt{2}+\sqrt{10}}{2}}; \; \cbox{yellow}{\si_1}, \cbox{yellow}{\si_1}]\approx$ &
$[2; \; \cbox{yellow}{\sqrt{2}}, \; \cbox{yellow}{\sqrt{10}}; \; \cbox{yellow}{\si_3},\;  \cbox{yellow}{\si_4} ]\approx$
	\smallskip	\\
 $[2; \cbox{yellow}{2.581}, \; \cbox{yellow}{2.581}; \; \cbox{yellow}{0.016}, \; \cbox{yellow}{0.016} ]$ &     
  $[2; \; \cbox{yellow}{1.414}, \; \cbox{yellow}{3.162}; \; \cbox{yellow}{0.041}, \; \cbox{yellow}{0.072} ]$  
\medskip \\  
\hline

$2\times [\sqrt{10}; \cbox{yellow}{\frac{2+\sqrt{2}}{2}}, \frac{\sqrt{2}+4}{2}; \cbox{yellow}{\si_1}, \cbox{yellow}{\si_2}]\approx$ &
$2\times [\sqrt{10}; \cbox{yellow}{\frac{2+\sqrt{10}}{2}},\frac{\sqrt{2}+4}{2}; \cbox{yellow}{\si_2}, \cbox{yellow}{\si_3}] \approx$
	\smallskip	\\
$[3.162; \cbox{yellow}{1.707}, 2.707; \cbox{yellow}{0.016}, \cbox{yellow}{0.029} ]$ &     
  $[3.162; \cbox{yellow}{2.581}, 2.707; \cbox{yellow}{0.029},  \cbox{yellow}{0.042} ]$  
\bigskip \\  
\hline
    
$[4; \; \frac{\sqrt{2}+\sqrt{10}}{2}, \; \frac{\sqrt{2}+\sqrt{10}}{2}; \; \si_2, \; \si_2] \approx$ &
$[4; \; \frac{\sqrt{2}+\sqrt{10}}{2}, \; \frac{\sqrt{2}+\sqrt{10}}{2}; \; \si_2, \; \si_2] \approx$
	\smallskip	\\
$[4; \; 2.288, \; 2.288; \; 0.029, \; 0.029 ]$ &     
$[4; \; 2.288, \; 2.288; \; 0.029, \; 0.029 ]$  
  \end{tabular}
  \caption{The Average Simplexwise Distributions from Definition~\ref{dfn:ASV+ASD+ASM} for the clouds $T,K\subset\R^2$ in Fig.~\ref{fig:cloud_T_strengths}, \ref{fig:cloud_K_strengths} consist of six vectors in $\R^5$. 
The coordinate-wise averages of six vectors in $\ASD$ give
the Average Simplexwise Moments $\ASM\in\R^5$ in Table~\ref{tab:ASM_TK}.}
  \label{tab:ASD_TK}
\end{table*}

The kite $K$ has two subsets subset $A=\{p_1,p_2\}$ and $A=\{p_2,p_3\}$ with $\SPD(A)$ consisting of one distance $\sqrt{2}$.
The average distance from $p_4$ and $p_2$ (or $p_3$) to $p_1$ is $\dfrac{4+\sqrt{10}}{2}$, which appeared for $T$ above.
However, the average distance from $p_4$ and $p_3$ to $p_2$ is $\dfrac{2+\sqrt{2}}{2}$, so its difference with $\dfrac{2+\sqrt{10}}{2}$ is \hl{highlighted} in the first row of Table~\ref{tab:ASD_TK}.
Also, the normalized strengths of triangles $\triangle p_2 p_3 p_4$ and $\triangle p_1 p_2 p_3$ in $K$ are $\si_2<\si_1$, which differ from the strengths $\si_0<\si_2$ for the similar subset $A$ in the trapezoid $T$.
\smallskip

Rows 2 and 3 in Table~\ref{tab:ASD_TK} compare $\ASV(T;A)$ and $\ASV(K;A)$ for 2-point subsets $A$ with inter-point distances $2$ (one subset in both $T,K$) and $\sqrt{10}$ (two subsets in each of $T,K$).
Row 4 in Table~\ref{tab:ASD_TK} contains the same vectors $\ASV(T;A)=\ASV(K;A)$ for $A$ with inter-point distance $4$, because each of $T,K$ has two equal triangles with the largest side $4$.
\smallskip
 
Table~\ref{tab:ASM_TK} shows the five coordinates of the Oriented Distance Moments $\ASM(T;1)$ for $t=1$ in Definition~\ref{dfn:ASV+ASD+ASM} obtained by coordinate-wise averaging six vectors in Table~\ref{tab:ASD_TK}.
In the first row both averages are equal to 
$\ASM_1=\dfrac{1}{6}( 2\sqrt{2}+2 + 2\sqrt{10}+4)=\dfrac{3+\sqrt{2}+\sqrt{10}}{3}\approx 2.53$.
\smallskip 

\begin{table*}[h!]
\centering
\begin{tabular}{l|l}
$\cbox{yellow}{T}$: Average Simplexwise Moment & 
$\cbox{yellow}{K}$: Average Simplexwise Moment \\
 \hline
  	$\ASM_1=\dfrac{3+\sqrt{2}+\sqrt{10}}{3}\approx 2.53$ & 
   	$\ASM_1=\dfrac{3+\sqrt{2}+\sqrt{10}}{3}\approx 2.53$ \\

   	$\ASM_2=\dfrac{2+\sqrt{2}+\sqrt{10}}{3}\approx \cbox{yellow}{2.19}$ & 
   	$\ASM_2=\dfrac{8+5\sqrt{2}+3\sqrt{10}}{12}\approx \cbox{yellow}{2.05}$ \\

   	$\ASM_3=\dfrac{4+\sqrt{2}+\sqrt{10}}{3}\approx \cbox{yellow}{2.86}$ & 
   	$\ASM_3=\dfrac{16+3\sqrt{2}+5\sqrt{10}}{12}\approx \cbox{yellow}{3.01}$ \\

   	$\ASM_4=\dfrac{5\si_1+\si_2}{6}\approx \cbox{yellow}{0.018}$ & 
   	$\ASM_4=\dfrac{5\si_2+\si_3}{6}\approx \cbox{yellow}{0.031}$ \\

   	$\ASM_5=\dfrac{\si_1+5\si_2}{6}\approx \cbox{yellow}{0.027}$ & 
   	$\ASM_5=\dfrac{\si_2+4\si_3+\si_4}{6}\approx \cbox{yellow}{0.045}$  \\

  \end{tabular}
  \caption{
Average Simplexwise Moments from Definition~\ref{dfn:ASV+ASD+ASM} for $T,K\subset\R^2$ in Fig.~\ref{fig:cloud_T_strengths}, \ref{fig:cloud_K_strengths} are obtained by averaging the six vectors in the Average Simplexwise Distributions $\ASD$ from Table~\ref{tab:ASD_TK}.
}
\label{tab:ASM_TK}
\end{table*}

The distance $L_\infty$ between the vectors $\ASM(T;1)$ and $\ASM(K;1)$ in $\R^5$ is $\dfrac{\si_2+4\si_3+\si_4}{6}-\dfrac{\si_1+5\si_2}{6}=\dfrac{\si_4+4\si_3-4\si_2-\si_1}{6} \approx 0.018$.
\bs
\end{exa}

The vector $\ASD(C;A)$ in Definition~\ref{dfn:ASV+ASD+ASM} can be strengthened by including normalized strengths with signs, which  distinguish mirror images.
In this case, $T$ and $K$ would be distinguished also in last row 4 of Table~\ref{tab:ASD_TK} by the pairs of triangles with the common edge of length 4, because $T$ has these triangles on one side of the edge, but $K$ has them on different sides of the edge.
However, including signs creates a dependence on ordering points in a subset $A$ of a cloud $C$, which makes the invariant larger.
\smallskip

The rest of this section will prove lower bounds for $\EMD_\infty$ on $\OSD$s and $\SCD$s in terms of simpler invariants from Definition~\ref{dfn:ASV+ASD+ASM} and its centered version in Definition~\ref{dfn:ACV+ACD+ACM} later.
Recall that the matrix $D(A)$ of pairwise distances in Definition~\ref{dfn:D(A)+M(C;A)} and vector $\SPD(A)$ of sorted pairwise distances in Definition~\ref{dfn:ASV+ASD+ASM} make sense for a sequence $A$ of points in any metric space. 

\begin{prop}[$L_\infty$ bounds for sorted vectors]
\label{prop:Linf_bounds}
\textbf{(a)}
For any vector $\vb{v}=(v_1,\dots,v_m)$, write its coordinates in increasing order $v_{\be(1)}\leq\dots\leq v_{\be(m)}$ for a permutation $\be$ of $1,\dots,m$.
Set $\sort(\vb{v})=(v_{\be(1)},\dots, v_{\be(m)})$.
Then $L_\infty(\vb{u},\vb{v})\geq L_\infty(\sort(\vb{u}),\sort(\vb{v}))$ for any vectors $\vb{u},\vb{v}\in\R^m$. 
\medskip

\noindent
\textbf{(b)}
Any sequences $A,A'$ of $n$ ordered points in a metric space satisfy the inequality $L_\infty\big(D(A),D(A')\big)\geq L_\infty\big(\SPD(A),\SPD(A')\big)$.
\end{prop}
\begin{proof}
\textbf{(a)}
Set $\ep=L_\infty(\vb{u},\vb{v})$.
For each $k=1,\dots,m$, we will prove that the $k$-th sorted coordinates $\sort_k(\vb{u})$ and $\sort_k(\vb{v})$ are $\ep$-close.
Assume by contradiction that $\sort_k(\vb{u})<\sort_k(\vb{v})-\ep$.
Then $\sort(\vb{u})$ has at least $k$ coordinates less than $\sort_k(\vb{v})-\ep$.
Hence all their $\ep$-close coordinates in $\vb{v}$ are less than $\sort_k(\vb{v})$, so the $k$-th ordered coordinate in $\sort(\vb{u})$ is less than $\sort_k(\vb{v})$, which contradicts the choice of $\sort_k(\vb{v})$.
The opposite assumption $\sort_k(\vb{u})>\sort_k(\vb{v})+\ep$ similarly leads to a contradiction.
Then $|\sort_k(\vb{u})-\sort_k(\vb{v})|\leq\ep$ for all $k=1,\dots,m$ implies that  $L_\infty(\vb{u},\vb{v})\leq\ep$.
\medskip

\noindent
\textbf{(b)}
Flattening the upper triangle of $D(A)$ to a vector $\vb{v}$ of $m=\dfrac{n(n-1)}{2}$ indexed distances gives $\sort(\vb{v})=\SPD(A)$.
It remains to use part (a).
\end{proof}

\begin{prop}[$\EMD$ respects bounds of base metrics]
\label{prop:EMD_respects_base_metrics}
In the notations of Definition~\ref{dfn:EMD}, let distributions $B=\{B_1,\dots,B_k\}$ and $C=\{C_1,\dots,C_l\}$ have base metrics $d,d'$ such that $d(B_i,C_j)\geq d'(B_i,C_j)$ for all $i,j$.
Then the Earth Mover's Distances $\EMD,\EMD'$ based on the base metrics $d,d'$, respectively, satisfy $\EMD(B,C)\geq\EMD'(B,C)$.
\end{prop}
\begin{proof}
Let $f_{ij}\in[0,1]$ be optimal flows that minimize the Earth Mover's Distance $\EMD(B,C)=\sum\limits_{i,j} f_{ij} d(B_i,C_j)$.
The given inequality $d\geq d'$ implies that
$\EMD(B,C)\geq \sum\limits_{i,j} f_{ij} d'(B_i,C_j)$.
Taking the minimum of the right hand side over all $f_{ij}\in[0,1]$ subject to the weight conditions in Definition~\ref{dfn:EMD} gives
$\EMD(B,C)\geq\EMD'(B,C)$.
\end{proof}

Recall that Definition~\ref{dfn:ASV+ASD+ASM} introduced the vectors
$\ar{M}(C;A)$ and $\ar{S}(C;A)$ of averages $\bar d(q)=\dfrac{1}{n}\sum\limits_{k=1}^n|\vb{q}-\vb{p_k}|$ and normalized strengths $\dfrac{\si(A\cup\{q\})}{\la_n}$, respectively, which are written in increasing order for all points $q\in C\setminus A$.

\begin{lem}[lower bounds for $\EMD_\infty$ on clouds]
\label{lem:EMD_bounds_clouds}
\textbf{(a)}
Let $C\subset\R^n$ be any cloud  of $m$ unordered points and a base sequence $A\subset C$ of $n$ ordered points $p_1,\dots,p_n$.
For any point $q\in C\setminus A$, let $\check d(q)\in\R^n$ be the point with coordinates $|\vb{q}-\vb{p_i}|$ for $i=1,\dots,n$.
Denote by $\check M(C;A)\subset\R^n$ the cloud of $m-n$ unordered points $\check d(q)$ for all $q\in C\setminus A$.
In the notations of Definition~\ref{dfn:ORD+OCD_BD+EMD}, any other cloud $C'\subset\R^n$ with a base sequence $A'\subset C'$ of $n$ ordered points satisfies the inequalies
$\EMD_\infty\big(\hat M(C;A),\hat M(C';A')\big)\geq$ 
$$\max\Big\{\EMD_\infty\big(\check M(C;A),\check M(C';A')\big),\EMD_\infty\big(\ar{S}(C;A),\ar{S}(C';A')\big)\Big\},$$

\noindent
$\EMD_\infty\big(\check M(C;A),\check M(C';A')\big)\geq 
\EMD_\infty\big(\ar{M}(C;A),\ar{M}(C';A')\big)$, and
$$\EMD_\infty\big(\ORD(C;A),\ORD(C';A')\big)\geq \EMD_\infty\big(\ASV(C;A),\ASV(C';A')\big).$$

\noindent
\textbf{(b)}
For any Minkowski metric $L_q$ on $\R^n$ for $q\in[1,+\infty]$ and clouds $C,C'\subset\R^n$, the Earth Mover's Distance $\EMD_q$ based on $L_q$ has the lower bound in terms of weighted centers of mass:
$\EMD_q(C,C')\geq L_q(O(C),O(C'))$.
\end{lem}
\begin{proof}
\textbf{(a)}
Arbitrarily order the points of $C\setminus A$ and $C'\setminus A'$ as $q_1,\dots,q_{m-n}$ and $q'_1,\dots,q'_{m'-n}$, respectively.
Then all points of the clouds $\hat M(C;A)$ and $\hat M(C';A')$ are ordered as $\hat d(q_1),\dots,q_{m-n}$ and $q'_1,\dots,q'_{m'-n}$, respectively.
All values of the sets $S(C;A)$ and $S(C';A')$ become ordered as 
$s_1(C;A),\dots,s_{m-n}(C;A)$ and $s_1(C';A'),\dots,s_{m'-n}(C';A')$, respectively.
\smallskip

Due to Proposition~\ref{prop:EMD_respects_base_metrics}(a), to prove the first lower bound, it suffices to check two similar inequalities for the base metric in $\R^{n+1}$, i.e.
$$L_\infty\big(\check d(q_i),\check d(q'_j)\big)\leq 
L_\infty\big(\hat d(q_i),\hat d(q'_j)\big)\geq
\big|s_i(C;A)-s_j(C';A')\big|.$$
The former inequality holds, because each $L_\infty(\hat d(q_i),\hat d(q'_j))$ is the maximum absolute difference between corresponding coordinates of the points $\hat d(q_i),\hat d(q'_j)\in\R^{n+1}$ whose first $n$ coordinates form $\check d(q_i),\check d(q'_j)\in\R^{n}$.
The latter inequality will follow below form the fact that the last $(n+1)$-st coordinates are normalized signed strengths $\dfrac{s(A\cup\{q_i\})}{\la_n}$ and $\dfrac{s(A'\cup\{q'_j\})}{\la_n}$, respectively.
By Definition~\ref{dfn:ORD+OCD_BD+EMD}, each $s(A\cup\{q_i\})=s(q_i;A)\si(A\cup\{q_i\})$ equals the strength up to the sign factor $s(q_i;A)$ in Definition~\ref{dfn:D(A)+M(C;A)}.
Then 
$$\begin{array}{l}
|s(A\cup\{q_i\})-s(A'\cup\{q'_j\})|\geq|\si(A\cup\{q_i\})-\si(A'\cup\{q'_j\})|, \text{ and hence}\\
L_\infty(\hat d(q_i),\hat d(q'_j))\geq 
\dfrac{|s(A\cup\{q_i\})-s(A'\cup\{q'_j\})|}{\la_n}\geq 
|s_i(C;A)-s_j(C';A')|.
\end{array}$$
The last expression is the distance $L_\infty$ between indexed elements of the unordered distributions $S(C;A)$ and $S(C';A')$ of normalized strengths.
Since $\EMD_\infty$ is independent of an order of objects, these distributions can be replaced with the vectors $\ar{S}(C;A)$ and $\ar{S}(C';A')$ from Definition~\ref{dfn:ASV+ASD+ASM}. 
\medskip

Now we will prove the lower bound involving the vectors $\ar{M}(C;A)$ and $\ar{M}(C';A')$ of average distances.
Each point $\check d(q_i)\in\R^n$ consists of $n$ coordinates $\check d_k(q_i)=|\vb{q_i}-\vb{p_k}|$ for $k=1,\dots,n$.
The base metric equals Chebyshev distance $L_\infty\big(\check d(q_i),\check d(q_j)\big)$ has the lower bound equal to the average of absolute differences between corresponding coordinates, i.e.
$$L_\infty\big(\check d(q_i),\check d(q_j)\big)
\geq \dfrac{1}{n}\sum\limits_{k=1}^n | d_k(q_i)- d_k(q'_j)|\geq
\dfrac{1}{n}\left|\sum\limits_{k=1}^n d_k(q_i)-\sum\limits_{k=1}^n d_k(q'_j)\right|=$$
$=|\bar d(q_i)-\bar d(q'_j)|$, which is the distance $L_\infty$ between corresponding coordinates of $\ar{M}(C;A)$ and $\ar{M}(C';A')$, whose re-ordering does not affect the Earth Mover's Distance.
The proved inequality between the base metrics gives rise to the same inequality between $\EMD$s by Proposition~\ref{prop:EMD_respects_base_metrics}(a). 
\medskip

To prove the final inequality, 
let $\xi\in S_n$ be the permutation that minimizes $\EMD_\infty[\xi]$.
By re-ordering the base sequence $A\subseteq C$, we can assume that $\xi$ is trivial.
Then $\EMD_\infty\big(\ORD(C;A),\ORD(C';A')\big)$ equals
$$e=\max\big\{ L_{\infty}\big(D(A),D(A')\big), \EMD_\infty\big(\hat M(C;A),\hat M(C';A')\big) \big\}.$$
The first lower bound $e\geq L_{\infty}\big(D(A),D(A')\big)
\geq L_\infty\big(\SPD(A),\SPD(A')\big)$ holds by 
Proposition~\ref{prop:Linf_bounds}(b).
The other lower bounds 
were proved above:
$$\begin{array}{l}
e\geq \EMD_\infty\big(\hat M(C;A),\hat M(C';A')\big)\geq \\
\max\Big\{\EMD_\infty\big(\ar{M}(C;A),\ar{M}(C';A')\big),\EMD_\infty\big(\ar{S}(C;A),\ar{S}(C';A')\big)\Big\}.
\end{array}$$

\noindent
\textbf{(b)}
It was proved in \cite[Theorem~3]{cohen1997earth}.
\end{proof}

In Lemma~\ref{lem:EMD_bounds_clouds}(a), if given clouds $C,C'\subset\R^n$ consists of the same number $m$ of points, then Proposition~\ref{prop:EMD_sets1D} provides the explicit formulas: 
$$\begin{array}{l}
\EMD_\infty\big(\ar{M}(C;A), \ar{M}(C';A')\big)=\dfrac{1}{m}L_1(\ar{M}(C;A), \ar{M}(C';A')\big), \\ \\
\EMD_\infty\big(\ar{S}(C;A), \ar{S}(C';A')\big)=\dfrac{1}{m}L_1(\ar{S}(C;A), \ar{S}(C';A')\big).
\end{array}$$

\begin{thm}[properties of metrics on $\ASD$s and $\ASM$s] 
\label{thm:ASD+ASM_metrics}
\textbf{(a)}
For any $t\geq 1$ and a cloud $C\subset\R^n$ of $m$ unordered points, the
Average Simplexwise Distribution $\ASD(C)$ and the Average Simplexwise Moment $\ASM(C;t)$ in Definition~\ref{dfn:ASV+ASD+ASM} can be computed in time $O(n^3m^{n+1}\log m)$.
\medskip

\noindent
\textbf{(b)}
For any clouds $C,C'\subset\R^n$ of up to $m$ points, the Earth Mover's Distance between $\ASD(C),\ASD(C')$s of a maximum size $l\leq\dbinom{m}{n}=O(m^n)$ in Definition~\ref{dfn:ASV+ASD+ASM}(b) can be computed in time $O(l^{3}\log l)$ and provides the lower bound
$\EMD_\infty\big(\OSD(C),\OSD(C')\big)\geq\EMD_\infty\big(\ASD(C),\ASD(C')\big)$, hence $\ASD(C)$ has Lipschitz constant $2$ under any perturbations  of $C$.
\medskip

\noindent
\textbf{(c)}
For any $t\geq 1$ and clouds $C,C'\subset\R^n$ with up to $m$ points and pre-computed $\ASM$s, 
$\EMD_\infty(\ASM(C;t);\ASM(C';t))$ in Definition~\ref{dfn:ASV+ASD+ASM}(c) is computed in extra time $O((n^6+m^3)\log m)$ and provides the lower bound 
$$\EMD_\infty\big(\ASD(C),\ASD(C')\big)\geq \EMD_\infty\big(\ASM(C;1),\ASM(C';1)\big),$$
hence $\ASM(C)$ has Lipschitz constant $2$ under any perturbations  of $C$.
\end{thm}
\begin{proof}
\textbf{(a)}
To compute $\OSD(C)$, we consider $k=\dbinom{m}{n}=\dfrac{m!}{n!(m-n)!}=O(m^n)$ unordered subsets $A\subseteq C$ of $n$ points.
The vector $\SPD(A)$ of $\dfrac{n(n-1)}{2}$ sorted pairwise distances can be computed in time $O(n^2\log n)$.
The vector $\ar{M}(C;A)$ consists of $m-n$ sorted averages based on $n(m-n)$ distances, all computable in time $O(mn+m\log m)$.
The vector $\ar{S}(C;A)$ consists of $m-n$ average strengths, each computable in time $O(n^3)$ by Theorem~\ref{thm:strength_properties}. 
Then $\ASV(C;A)=[\SPD(A);\ar{M}(C;A);\ar{S}(C;A)]$ is computed in time $O(n^3m\log m)$.
Multiplying the last time by the number $k=O(m^n)$ of $n$-point subsets $A\subseteq C$, the time of $\ASD(C)$ becomes $O(n^3m^{n+1}\log m)$.
After that, the $t$-th moment $\ASM(C;t)$ needs extra time $O(m^n)$ for each of $\dfrac{n(n-1)}{2}+2(m-n)$ coordinates, which keeps the total time the same.
\medskip

\noindent
\textbf{(b)}
In the notations of Definition~\ref{dfn:ASV+ASD+ASM}, the distance $L_\infty(\SPD(A),\SPD(A))$ on $\dfrac{n(n-1)}{2}$ ordered distances needs time $O(n^2\log n)$.
The vectors $\ar{M}(C;A)$, $\ar{M}(C';A')$, $\ar{S}(C;A)$, $\ar{S}(C';A')$ have $O(m)$ coordinates.
Hence any distance $\EMD_\infty$ between them needs time $O(m^3\log m)$, see \cite{fredman1987fibonacci,goldberg1987solving}.
Since $n\leq m$, the time for $\EMD_\infty\big(\ASV(C;A),\ASV(C;A)\big)$ is $O(m^3\log m)$.
\smallskip

If the distributions $\ASD(C),\ASD(C')$ consist of $l\leq\dbinom{m}{n}=O(m^n)$ unordered $\ASV$s based on subsets of $n$ points, $\EMD_\infty\big(\ASD(C),\ASD(C')\big)$ can be computed in time $O(l^3\log l)=O(nm^{3n}\log m)$.
Both $\OSD(C)$ and $\ASD(C)$ can be considered unordered distributions of equally weighted pairs $(C;A)$ with the base metrics
$d=\EMD_\infty\big(\ORD(C;A),\ORD(C';A')\big)$, $d'=\EMD_\infty\big(\ASV(C;A),\ASV(C';A')\big)$, respectively.
Then $d\geq d'$ from Lemma~\ref{lem:EMD_bounds_clouds}(a) implies the same inequality for $\EMD_\infty$ by Proposition~\ref{prop:EMD_respects_base_metrics}. 
This inequality and continuity of $\EMD_\infty$ on $\OSD$s with Lipschitz constant $2$ in Theorem~\ref{thm:OSD+SCD_continuous} imply the same continuity of $\EMD_\infty$ on $\ASD$s.  
\medskip

\noindent
\textbf{(c)}
For a fixed $t\geq 1$, the Average Simplexwise Moment $\ASM(C;t)$ consist of the vector $\mu_t\{\SPD(A)\}$ of $\dfrac{n(n-1)}{2}=O(n^2)$ coordinates and the vectors $\mu_t\{\ar{M}(C;A)\}$ and $\mu_t\{\ar{S}(C;A)\}$ of $m-n=O(m)$ coordinates each.
Computing $\EMD_{\infty}$ between vectors of $k$ coordinates needs time $O(k^3\log k)$, which gives the overall time $O((n^6+m^3)\log m)$.
\smallskip
 
To prove the lower bound for $\EMD_\infty$ on $\ASM(C;1)$, separately consider the distributions $\SPD(A)$, $\ar{M}(C;A)$, and $\ar{S}(C;A)$ of subvectors.
Then apply Lemma~\ref{lem:EMD_bounds_clouds}(b) to the centers of each distribution as follows:
$$\begin{array}{l}
\EMD_\infty\big(\SPD(A),\SPD(A')\big)\geq L_\infty \big(\mu_1\{\SPD(A)\},\mu_1\{\SPD(A')\}\big), \\
\EMD_\infty\big(\ar{M}(C;A),\ar{M}(C';A')\big)\geq L_\infty \big(\mu_1\{\ar{M}(C;A)\},\mu_1\{\ar{M}(C';A')\}\big), \\
\EMD_\infty\big(\ar{S}(C;A),\ar{S}(C';A')\big)\geq L_\infty \big(\mu_1\{\ar{S}(C;A)\},\mu_1\{\ar{S}(C';A')\}\big).
\end{array}$$ 
The inequalities imply the similar required inequality for the overall distances $\EMD_\infty$, which were introduced in Definition~\ref{dfn:ASV+ASD+ASM}(b,c) as  the maxima of  three $\EMD_\infty$ distances for both $\ASD$ and $\ASM$.
The proved lower bounds and the continuity of $\EMD_\infty$ on $\OSD$s with Lipschitz constant $2$ in Theorem~\ref{thm:OSD+SCD_continuous} imply the same continuity of $\EMD_\infty$ on $\ASM$s.  
\end{proof}

Though the rough time estimate $O(nm^{3n}\log m)$ for $\EMD_\infty$ on $\ASD$s in Theorem~\ref{thm:ASD+ASM_metrics}(b) is asymptotically the same as for $\EMD_\infty$ on $\OSD$s in Theorem~\ref{thm:OSD_metrics}(b), averaging over points of $C\setminus A$ likely makes $\ASD(C)$ smaller than $\OSD(C)$ for symmetric clouds, as in Table~\ref{tab:ASD_TK}.
\smallskip

The time $O((n^6+m^3)\log m)$ for $\EMD_\infty$ on vectors $\ASM$s
is certainly smaller than $O(nm^{3n}\log m)$ for practical dimensions $n=2,3$.
We reduce computational times using the center of mass as an extra base point.

\begin{table}[h!]
\centering
\begin{tabular}{lll}
$\ASD(C)$ & Average Simplexwise Distribution  & Definition~\ref{dfn:ASV+ASD+ASM} \\
$\ASM(C;t)$ & Average Simplexwise $t$-th Moment  & Definition~\ref{dfn:ASV+ASD+ASM} \\x
$\ASV(C;A)$ & Average Simplexwise Vector for $A\subset C$  & Definition~\ref{dfn:ASV+ASD+ASM} \\
$\ACD(C)$ & Average Centered Distribution  & Definition~\ref{dfn:ACV+ACD+ACM} \\
$\ACM(C;t)$ & Average Centered $t$-th Moment  & Definition~\ref{dfn:ACV+ACD+ACM} \\
$\ACV(C;A)$ & Average Centered Vector for $A\subset C$  & Definition~\ref{dfn:ACV+ACD+ACM} \\
$\BMD(\Ga)$ & Bottleneck Matching Distance of a graph  & Definition~\ref{dfn:BMD} \\
$\EMD(B,C)$ & Earth Mover's Distance on distributions  & Definition~\ref{dfn:EMD} 
\end{tabular}
\caption{Acronyms and references for all concepts in sections~\ref{sec:OSD+SCD_metrics} and \ref{sec:OSD+SCD_moments}.}
\label{tab:acronyms_SCD_metrics}
\end{table} 

\begin{dfn}[Average Centered Moments $\ACM(C;t)$]
\label{dfn:ACV+ACD+ACM}
Let $C\subset\R^n$ be a cloud of $m$ unordered points with a center of mass $O(C)$.
For a base sequence $A\subset C$ of $n-1$ ordered points $p_1,\dots,p_{n-1}$, let the extended sequence $\ti A=A\cup O(C)$ include $O(C)$ as the last $n$-th point.
The vector $\ar{D}(C;A)$ is formed by $n-1$ distances $|\vb{p_i}-\vb{O(C)}|$, $i=1,\dots,n-1$, written in increasing order.
The vector $\ar{O}(C;A)$ is formed by $m-n+1$ distances $|\vb{q}-\vb{O(C)}|$,  written in increasing order, for all points $q\in C\setminus A$.
\medskip

\noindent
\textbf{(a)}
In the notations of Definition~\ref{dfn:ASV+ASD+ASM}, the \emph{Average Centered Vector} 
$$\ACV(C;A)=[\SPD(A);\ar{D}(C;A);\ar{M}(C;A);\ar{O}(C;A);\ar{S}(C;\ti A)]$$
has $\dfrac{n(n-1)}{2}+3(m-n+1)$ coordinates.
For any cloud $C'\subset\R^n$ of $m'$ unordered points 
with a base sequence $A'$ of $n-1$ points, 
define $\EMD_\infty\big(\ACV(C;A),\ACV(C';A')\big)$ as the maximum of the five distances
\smallskip

\hspace*{-8mm}
$\begin{array}{l}
L_\infty\big(\SPD(A),\SPD(A')\big), \\
L_\infty\big(\ar{D}(C;A),\ar{D}(C';A')\big), \\
\EMD_\infty\big(\ar{M}(C;A), \ar{M}(C';A')\big), \\
\EMD_\infty\big(\ar{O}(C;A), \ar{O}(C';A')\big), \\
\EMD_\infty\big(\ar{S}(C;\ti A), \ar{S}(C';\ti A')\big).
\end{array}$
\medskip

\noindent
\textbf{(b)}
The \emph{Average Centered Distribution} $\ACD(C)$ is the unordered set of equally weighted vectors $\ACV(C;A)$ for all $\dbinom{m}{n-1}$  unordered subsets $A\subseteq C$ of $n-1$ points.
Use the base metric $\EMD_\infty$ on $\ACV$s in part~(a) to define
the Earth Mover's Distance $\EMD_\infty\big(\ACD(C),\ACD(C')\big)$. 
\medskip

\noindent
\textbf{(c)}
For $t\geq 1$, the \emph{Average Centered Moment} $\ACM(C;t)=\mu_t\{\ACD(C)\}$ is the $t$-th moment of $\ACD(C)$ of equally weighted vectors $\ACV(C;A)$ over all $\dbinom{m}{n-1}$ unordered subsets $A\subseteq C$ of $n-1$ points, i.e.
$\ACM(C;t)=$
$$\big[\mu_t\{\SPD(A)\};\mu_t\{\ar{D}(C;A)\};\mu_t\{\ar{M}(C;A)\}; \mu_t\{\ar{O}(C;A)\};\mu_t\{\ar{S}(C;\ti A)\}\big].$$
For any clouds $C,C'\subset\R^n$, define $\EMD_\infty\big(\ACM(C;t),\ACM(C';t)\big)$ as the maximum of the five metrics $\EMD_\infty$ between constituent subvectors.
\bs
\end{dfn}

$\ACV(C;A)$ in Definition~\ref{dfn:ACV+ACD+ACM}(a) uses the extended base sequence $\ti A$ only in the last vector $\ar{S}(C;A)$ because the strengths need $n+1$ points from $\ti A\cup\{q\}$.
Example~\ref{exa:invariants_regular_clouds} simplifies the complete invariant $\SCD$ of regular clouds from Example~\ref{exa:SCD} to distributions and moments from Definition~\ref{dfn:ACV+ACD+ACM}.

\begin{exa}[centered distributions for regular clouds $C_3,C_4\subset\R^2$]
\label{exa:invariants_regular_clouds}
\textbf{(a)}
The cloud $C_3\subset\R^2$ at the vertices $p_1,p_2,p_3$ of an equilateral triangle has the center of mass $O(C_3)=\vb*{0}$ at the origin.
For any base sequence $A_i=\{p_i\}$ of one point, and another point $p_j\neq p_i$, the triangle $\triangle p_i p_j \vb*{0}$ has $\Vol=\dfrac{\sqrt{3}}{4}$, half-perimeter $p=\dfrac{2+\sqrt{3}}{2}$, and the normalized strength $\bar \si_3=\dfrac{\Vol^2}{p^3\la_2 }=\dfrac{(\sqrt{3})^2 2^3}{4^2(2+\sqrt{3})^3\sqrt{3}}=\dfrac{\sqrt{3}}{2(2+\sqrt{3})^3}\approx 0.017$.
Example~\ref{exa:SCD}(a) has the Oriented Centered Distribution $\OCD(C_3;A_i)=\left[1,\left( \begin{array}{cc} 
\sqrt{3} & \sqrt{3} \\
1 & 1  \\ 
-1 & +1
\end{array}\right) \right]$.
By Definition~\ref{dfn:ACV+ACD+ACM}, $\SPD(A_i)$ is empty, $\ar{D}(C_3;A_i)$ consists of one distance $|p_i|=1$, the vector $\ar{M}(C_3;A_i)$ consists of two distances $|\vb{q}-\vb{p_i}|=\sqrt{3}$, the vector $\ar{O}(C_3;A_i)$ consists of two distances $|\vb{q}|=1$, and $\ar{S}(C_3;\ti A_i)$ consists of two normalized strengths $\dfrac{\si(\ti A_i\cup\{q\})}{\la_2}=\bar\si_3$ for both points $q\in C_3\setminus\{p_i\}$.
The Average Centered Vector is $\ACV(C_3;A_i)=[1;\sqrt{3},\sqrt{3};1,1;\bar\si_3,\bar\si_3]\in\R^7$.
Since $\ACV(C_3;A_i)$ is independent of $i=1,2,3$, the Average Centered Distribution $\ACD(C_3)$ and the Average Centered Moment $\ACM(C_3;1)$ coincide with $[1;\sqrt{3},\sqrt{3};1,1;\bar\si_3,\bar\si_3]$, and $\ACM(C_3;t)=\vb{0}$ for $t>1$.
\medskip

\noindent
\textbf{(b)}
The cloud $C_4\subset\R^2$ at the vertices $p_1,p_2,p_3,p_4$ of a square has the center of mass $O(C_4)=\vb*{0}$ at the origin.
For any base sequence $A_i=\{p_i\}$ of one point, and another point $p_{i\pm 1}$, where $i\pm 1$ is taken modulo 4 in $\{1,2,3,4\}$, the triangle $\triangle p_i p_{i\pm 1} \vb*{0}$ has $\Vol=\dfrac{1}{2}$, half-perimeter $p=\dfrac{2+\sqrt{2}}{2}$, and the normalized strength $\bar \si_4=\dfrac{\Vol^2}{p^3\la_2 }=\dfrac{2^3}{2^2(2+\sqrt{2})^3\sqrt{3}}=\dfrac{2}{\sqrt{3}(2+\sqrt{3})^3}\approx 0.022$.
Any other triangle $\triangle p_i p_{i+2} \vb*{0}$ is degenerate, so its strength is 0.
Example~\ref{exa:SCD}(b) has the Oriented Centered Distribution
$\OCD(C_4;\ti A_i)=\left[ 1, \left( \begin{array}{ccc} 
\sqrt{2} & 2 & \sqrt{2} \\
1 & 1 & 1 \\ 
-1 & 0 & +1
\end{array}\right) \right]$.
Then $\SPD(A_i)$ is empty, $\ar{D}(C_4;A_i)$ consists of one distance $|p_i|=1$, the vector $\ar{M}(C_4;A_i)$ consists of three distances $\sqrt{2},\sqrt{2},2$, the vector $\ar{O}(C_3;A_i)$ consists of three distances $|\vb{q}|=1$, and $\ar{S}(C_3;\ti A_i)$ consists of $0$ and two normalized strengths $\dfrac{\si(\ti A_i\cup\{q\})}{\la_2}=\bar\si_4$.
The Average Centered Vector is $\ACV(C_4;A_i)=[1;\sqrt{2},\sqrt{2},2;1,1,1;0,\bar\si_4,\bar\si_4]\in\R^{10}$.
The vector $\ACV(C_4;A_i)$ is independent of $i$, then $\ACD(C_4)$ and $\ACM(C_4;1)$ coincide with this vector. 
\medskip

\noindent
\textbf{(c)}
To compare the triangular and square clouds $C_3,C_4$, we compute the metric $\EMD_\infty\big(\ACD(C_3),\ACD(C_4)\big)$ as the maximum of three distances $\EMD_\infty$, which can be found by Proposition~\ref{prop:EMD_sets1D}:
$\EMD_\infty([1,1],[1,1,1])=0$,
$$\begin{array}{l}
\EMD_\infty([\sqrt{3},\sqrt{3}],[\sqrt{2},\sqrt{2},2])=
\dfrac{1}{3}(2-\sqrt{3})+\dfrac{2}{3}(\sqrt{3}-\sqrt{2})
\approx 0.3, \\
\EMD_\infty([\bar\si_{3},\bar \si_{3}],[0,\bar \si_{4},\bar \si_{4}])=
\dfrac{1}{3}(\bar\si_{3}-0)+\dfrac{2}{3}(\bar\si_{4}-\bar\si_{3})\approx 0.003,
\text{ so}
\end{array}$$
$\EMD_\infty\big(\ACD(C_3),\ACD(C_4)\big)=\EMD_\infty\big(\ACM(C_3;1),\ACM(C_4;1)\big)\approx 0.3$, which slightly differs from another $\EMD=\dfrac{5}{16}$ in \cite[Example~6.7]{kurlin2024polynomial}.
\bs
\end{exa}

\begin{thm}[properties of metrics on $\ACD$s and $\ACM$s] 
\label{thm:ACD+ACM_metrics}
\textbf{(a)}
For any $t\geq 1$ and a cloud $C\subset\R^n$ of $m$ unordered points, the
Average Centered Distribution $\ACD(C)$ and the Average Centered Moment $\ACM(C;t)$ in Definition~\ref{dfn:ACV+ACD+ACM} can be computed in time $O(n^3m^{n}\log m)$.
\medskip

\noindent
\textbf{(b)}
For any clouds $C,C'\subset\R^n$ of up to $m$ points, the Earth Mover's Distance between $\ACD(C),\ACD(C')$ of a maximum size $l\leq\dbinom{m}{n-1}=O(m^{n-1})$ in Definition~\ref{dfn:ACV+ACD+ACM}(b) can be computed in time $O(l^{3}\log l)$ and satisfies
$\EMD_\infty\big(\OCD(C),\OCD(C')\big)\geq\EMD_\infty\big(\ACD(C),\ACD(C')\big)$, hence $\ACD(C)$ has Lipschitz constant $2$ under any perturbations  of $C$.
\medskip

\noindent
\textbf{(c)}
For any $t\geq 1$ and clouds $C,C'\subset\R^n$ with up to $m$ points and pre-computed $\ACM$s, the distance 
$\EMD_\infty(\ACM(C;t);\ACM(C';t))$ in Definition~\ref{dfn:ACV+ACD+ACM}(c) can be computed in extra time $O((n^6+m^3)\log m)$ and satisfies
$$\EMD_\infty\big(\ACD(C),\ACD(C')\big)\geq \EMD_\infty\big(\ACM(C;1),\ACM(C';1)\big),$$
so $\ACM(C;1)$ has Lipschitz constant $2$ under any perturbations  of $C$.
\end{thm}
\begin{proof}
All arguments are similar to the proof of Theorem~\ref{thm:ASD+ASM_metrics}, where $n$ points in a base sequence $A$ are replaced with $n-1$ points, because the extra $n$-th point is the center of mass $O(C)$.
Lemma~\ref{lem:EMD_bounds_clouds}(a) is similarly updated with five constituent vectors in Definition~\ref{dfn:ACV+ACD+ACM}(a), instead of three.
\end{proof}

In conclusion, this paper solved Problem~\ref{pro:cloud_invariants_SCD} for all $n$-dimensional clouds $C\subset\R^n$ of $m$ unordered points by the complete invariant $\SCD$ in Theorem~\ref{thm:SCD_complete} with Lipschitz continuous metrics due to Theorems~\ref{thm:SCD_metrics} and \ref{thm:OSD+SCD_continuous}, all computable in polynomial times of $m$, for a fixed dimension $n$, with faster invariants and lower bounds for distances in Theorem~\ref{thm:ACD+ACM_metrics}.
\medskip

Problem~\ref{pro:cloud_invariants_SCD} has been extended with extra conditions in \cite[Problem~1.4.5]{anosova2025geometric} in the emerging area of Geometric Data Science, which develops complete and continuous invariants of real data beyond clouds \cite{widdowson2022average,kurlin2023simplexwise,anosova2024importance}.
\medskip

\noindent 
Acknowledgment.
The author thanks Daniel Widdowson for implementing the original invariants, see the short conference version \cite{widdowson2023recognizing} and the extended version \cite{kurlin2023strength} with early proofs of some results in this paper. 
 
\bibliographystyle{match} 
\bibliography{simplexwise-centered-distribution}

\end{document}